\documentclass[a4paper,11pt]{article}
\usepackage{amssymb}
\usepackage{amssymb}
\usepackage{amssymb}
\usepackage{graphicx}
\usepackage{amsfonts}
\usepackage{amsmath}
\usepackage{amssymb}
\usepackage{fancyhdr}
\usepackage{indentfirst}
\usepackage{titlesec}
\usepackage{hyperref}

\begin{document}

\newtheorem{Them}{Theorem}[section]
\newtheorem{Prop}{Proposition}[section]
\newtheorem{Remark}{Remark}[section]
\newtheorem{Lemma}{Lemma}[section]
\newenvironment{Proof}{\textbf{Proof:}}{\begin{flushright} $\Box$ \end{flushright}}
\newtheorem{Definition}{Definition}[section]
\newtheorem{Claim}{Claim}
\newtheorem{Cor}{Corollary}[section]



\title{\textbf{On the existence of min-max minimal torus}}
\author{Xin Zhou}
\date{}
\maketitle

\pdfbookmark[0]{On the existence of min-max minimal torus}{beg}


\begin{abstract}
In this paper, we will study the existence problem of minmax
minimal torus. We use classical conformal invariant geometric
variational methods. We prove a theorem about the existence of
minmax minimal torus in Theorem \ref{convergence theorem}.
Firstly we prove a strong uniformization result(Proposition \ref{uniformization for torus}) using method of \cite{AB}. Then we use this proposition to choose good parametrization for our
minmax sequences. We prove a compactification result(Lemma \ref{compactification}) similar to that of
Colding and Minicozzi \cite{CM}, and then give bubbling convergence results similar
to that of Ding, Li and Liu \cite{DLL}. In fact, we get an
approximating result similar to the classical deformation lemma(Theorem \ref{main theorem}).
\end{abstract}



\section{Introduction}

The existence problem of minimal surfaces is always an interesting
topic. We know the existence of minimizing minimal disk, i.e. the
classical Plateau problem (see Chapter 4 of \cite{CM1}) since 1931.
There are many results from that time. In general, a minimal surface
is a harmonic conformal branched immersion from a Riemann surface to
a compact Riemannian manifold. Most results only consider existence
of \textbf{area minimizing} minimal surfaces in a given homotopy
class. In particular, the existence of area-minimizing surfaces
has been proved for all genus in a suitable sense (cf. \cite{SU1}, \cite{SU2}, \cite{CT} etc.).\\

Besides minimizing minimal surfaces, we naturally ask whether there
exist \textbf{min-max} minimal surfaces. Here \textbf{min-max} means
the area of the minimal surfaces is just the min-max critical point
of the area functional in a homotopy class. In general, suppose $A$ is a
functional on a Banach manifold $\mathbf{M}$,
$\Omega=\big\{v(t):[0,1]\rightarrow \mathbf{M},\ v\in C^{0}([0, 1], \mathbf{M})\big\}$
the path space in $\mathbf{M}$ with $\sigma\in\Omega$. Then
$\mathcal{W}_{A}=\underset{\rho\in[\sigma]}{inf}\quad\underset{t\in[0,1]}{max}A\big(\rho(t)\big)$
is the min-max critical value in the homotopy class of $\rho$. It will be more complicated when
considering min-max minimal surfaces than the minimizing case. From
the point of view of variational method, the approximation sequences
will be one parameter families of mappings, which makes it difficult
to do compactification. J. Jost gave such an approach in his book
\cite{J1}. Recently Colding and Minicozzi \cite{CM} also gave such
an approach in the case of sphere using geometric variational
methods. They all used the bubble convergence of almost harmonic
mappings from closed surfaces given by Sacks and Uhlenbeck
\cite{SU1}. Colding and Minicozzi also found a good approximation
sequence which plays an important role in their proof of finite time
extinction of the Ricci flow.\\

We will extend Colding and Minicozzi' approach to the case of torus,
i.e. the existence of min-max minimal torus. In fact, we give a
stronger approximation for a special minimizing sequence.
\emph{Using notations in Section 2.1}, the main result is:

\begin{Them}\label{main theorem}
For any homotopically nontrivial path $\beta\in\Omega$, if $\mathcal{W}>0$, there exists
a sequence $(\rho_{n}, \tau_{n})\in[\beta]$, with
$\underset{t\in[0,1]}{max}E\big(\rho_{n}(t),
\tau_{n}(t)\big)\rightarrow\mathcal{W}$, and $\forall\epsilon>0$,
there exist $N$ and $\delta>0$ such that if $n>N$, then for any
$t\in(0, 1)$ satisfying:
\begin{equation}\label{energy condition}
E\big(\rho_{n}(t), \tau_{n}(t)\big)>\mathcal{W}-\delta,
\end{equation}

there are possibly a conformal harmonic torus
$u_{0}:T^{2}_{\tau_{0}}\rightarrow N$ and finitely many harmonic
sphere $u_{i}:S^{2}\rightarrow N$, such that:
\begin{equation}
d_{V}\big(\rho_{n}(t), \underset{i}{\cup}u_{i}\big)\leq\epsilon.
\end{equation}
\end{Them}

Here $d_{V}$ means varifold distance as in Appendix A in \cite{CM}.
It is a corollary of Theorem \ref{convergence theorem} and Appendix
A in \cite{CM}. It is a stronger approximation result than Theorem
1.14 of \cite{CM}. We use the energy condition inequality
\ref{energy condition} for the special sequence $\rho_{n}$, while
Theorem 1.14 of \cite{CM} use area condition.\\

In the case of torus, we have to include the variation of conformal
structures as discussed in \cite{SU2} and \cite{SY}. The analysis of
singularity in the bubble convergence will be more complicated than
in the case of sphere. We will give existence results similar to
that of Ding, Li and Liu \cite{DLL}. In the following, we will first
give our notations, and then give the sketch of this paper an the
end of Section \ref{sketch of the variational approach}.\\

\textbf{Acknowledgement.} This is part of my master degree thesis in
Peking University. I would like to sincerely thank my advisor Professor
Gang Tian for his longtime help and encouragement. I also would like
to thank Professor Weiyue Ding for several valuable talks with me about
this problem. I would like to thank Professor Richard Schoen, Professor
Tobias Colding and Professor William P. Minicozzi
\uppercase\expandafter{\romannumeral2} for their interest in this
work. I would like to thank Yalong Shi for carefully reading and suggestions on the paper.
Finally I would like to thank Professor Bin Xu for his patience to
listen to my primitive ideas of this paper during his visit to
Beijing International Center for Mathematical Research.


\section{Sketch of the variational methods for min-max minimal torus}

In the paper \cite{CM}, Colding and Minicozzi used variational
methods to give the existence of min-max minimal spheres. Let's
firstly sketch their idea. Let$(N, h)$ be the ambient space.
$\Omega=\bigg\{\gamma(t)\in C^{0} \Big([0,1],C^{0}\cap
W^{1,2}(\mathbf{S^{2}}, N)\Big)\bigg\}$ is the path space. Here for
all $\gamma(t)\in\Omega$, $\gamma(0)$, $\gamma(1)$ are constant
mappings. We call all such one parameter family of mappings
$\gamma(t)\in\Omega$ \textbf{paths} in the following. For
$\beta\in\Omega$, let $[\beta]$ be the homotopy class of $\beta$ in
$\Omega$. The min-max critical value is
$\mathcal{W}=\underset{\rho\in [\beta]}{inf}\quad\underset{t\in
[0,1]}{max}Area\big(\rho(t)\big)$. They want to learn the behavior
of critical points corresponding to $\mathcal{W}$. They firstly
chose an arbitrary minimizing sequence $\tilde{\gamma}_{n}(t)\in
[\beta]$, such that
$\lim_{n\rightarrow\infty}\underset{t\in[0,1]}{max}Area\big(\tilde{\gamma}_{n}(t)\big)=\mathcal{W}$.
Then they did almost conformal reparametrization for these paths to
get $\gamma_{n}(t)\in [\beta]$ which are almost conformal, i.e.
$E(\gamma_{n}(t))-Area(\gamma_{n}(t))\rightarrow 0$. Finally they
perturbed $\gamma_{n}(t)$ to $\rho_{n}(t)$ by local harmonic
replacement so that the new paths $\rho_{n}(t)$ have certain
compactness. The existence of min-max minimal spheres follows from
this construction and Sacks and Uhlenbeck's bubbling compactness
\cite{SU1}.\\

We want to extend the min-max variational method given by Colding
and Minicozzi to the case of torus $T^{2}$. The difference between
sphere and torus is that torus has more than one conformal
structures, while the conformal structure of sphere is unique.
Generally speaking, the pull-back metrics of the mappings on the
area minimizing sequence of paths will correspond to different
conformal structures. It is natural to include the variance of the
conformal structures in the min-max construction. In fact, we need
to consider the Teichm$\ddot{u}$ller space of torus in order to
maintain the homotopy class of the paths as discussed in \cite{SY}.
It will be difficult to do both conformal reparametrization and
compactification, and we must also consider whether the
corresponding conformal structures converge. Fortunately, the
Teichm$\ddot{u}$ller space of $T^{2}$ is easy to manipulate, and the
singularity arising from the absence of compactness of conformal
structures has been given in \cite{DLL} by Ding, Li and Liu.


\subsection{Teichm$\ddot{u}$ller space of torus and the notations}

We know that any flat torus $T^{2}$ can be viewed as the quotient
space of $\mathbb{C}$ moduled by a lattice generated by bases
$\{\omega_{1}, \omega_{2}\}$. After some conformal linear
transformation, we can assume $\omega_{1}=1$, and
$\omega_{2}=\tau=\frac{w_{2}}{w_{1}}$, where $\tau$ lies in the
upper half plane $\mathbb{H}$. In fact \emph{the
Teichm$\ddot{u}$ller spaces of torus $\mathcal{T}_{1}$}, is just the
upper half plane $\mathbb{H}$. We call each element
$\tau\in\mathcal{T}_{1}$ a \textbf{mark}, and denote $\tau$ by a
marked torus $(T^{2}, \tau)$ as in Definition 2.7.2 of \cite{J2},
which means a torus by gluing edges of the lattice $\{1, \tau\}$
with the plane metric $dzd\overline{z}$. Denoting
$\tau=\tau_{1}+\sqrt{-1}\tau_{2}$, we have another normalization
such that the area of the corresponding torus $Area(\{\omega_{1},
\omega_{2}\})=1$, i.e. by letting
$\omega_{1}=\frac{1}{\sqrt{\tau_{2}}}$,
$\omega_{2}=\frac{\tau_{1}}{\sqrt{\tau_{2}}}+\sqrt{-1}\sqrt{\tau_{2}}$.
Let $T^{2}_{0}$ be the marked torus $(T^{2}, \sqrt{-1})$, then there
is \emph{a natural diffeomorphism $i_{\tau}$} from $(T^{2}, \tau)$
to $(T^{2}, \sqrt{-1})$, which is the quotient map of the linear map
of $\mathbb{C}$ keeping $1$ and sending $\tau$ to $\sqrt{-1}$. So we
can also denote $\tau\in\mathcal{T}_{1}$ as $(T^{2}_{\tau},
i_{\tau})$ as in page 78 of \cite{J2}. We will show that every
metric on $T^{2}_{0}$ is conformal to a marked torus $(T^{2},
\tau)$, while keeping the conformal homeomorphism in the
homotopy class of $i_{\tau}^{-1}$.\\

\begin{Definition}
 Let $\tilde{\Omega}=\bigg\{ \big(\gamma(t), \tau(t)\big);
\gamma(t)\in C^{0}\Big([0,1], C^{0}\cap W^{1,2}\big( (T^{2},
\tau(t)) , N \big)\Big), \tau(t)\in C^{0}\big([0,1],
\mathcal{T}_{1}\big) \bigg\}$, and $\Omega=\bigg\{ \gamma(t)\in
C^{0}\Big([0,1], C^{0}\cap W^{1,2}(T^{2}_{0}, N) \Big) \bigg\}$.
We assume $\gamma(0), \gamma(1)$ are constant
mapping or map the torus to some circles in $N$. And $\tau(0),
\tau(1)=\sqrt{-1}$, if mappings on the endpoints are constant mappings, and not
restrained if not.
\end{Definition}

We use varying domains $(T^{2}, \tau(t))$ in the
definition of $\tilde{\Omega}$, and there are two ways to understand
this: we can pull back all $\gamma(t)$ to $T^{2}_{0}$ by
$i_{\tau_{t}}^{-1}$ and the continuity is defined w.r.t the same
domain $T^{2}_{0}$; we can also consider $\gamma(t)$ as defined on a
large ball of $\mathbb{C}$ containing all parallelograms generated
by $\{1, \tau(t)\}$, and continuity is defined w.r.t. the plane
ball. Since $\tau(t)$ is continuous, the two definitions are equivalent.
Here $\tilde{\Omega}$ and $\Omega$ are our variational spaces.\\

For the area functional, we only need to consider variational
problem in the space $\Omega$, since changing domain metrics will
not change the area. But for energy functional, different conformal
structures may lead to different energy, so we have to consider
variational problem in the space $\tilde{\Omega}$. Fix a
homotopically nontrial path $\beta\in\Omega$, $\big(\beta(t),
\tau_{0}(t) \big)\in\tilde{\Omega}$.\footnote{Here
$\tau_{0}(t)\equiv\sqrt{-1}$.} Let $[\beta]$ be the homotopy class
of $\beta$ in $\Omega$. Since path $\gamma(t)\in\tilde{\Omega}$ may
have different domains $T^{2}_{\tau(t)}$, \emph{the homotopy
equivalence $\alpha\sim\beta$ of
$\alpha(t):T^{2}_{\tau(t)}\rightarrow N$ and
$\beta(t):T^{2}_{\tau^{\prime}(t)}\rightarrow N$} is defined as
follows. We can identify $T^{2}_{\tau(t)}$, and
$T^{2}_{\tau^{\prime}(t)}$ to $T^{2}_{0}$ by $i_{\tau(t)}$ and
$i_{\tau^{\prime}(t)}$, then we can view $\alpha(t)$ and $\beta(t)$
as mappings defined on the same domain $T^{2}_{0}$ and hence define
their homotopy equivalence.

\begin{Definition}
Let $\mathcal{W}=\underset{\rho\in[\beta]}{inf}\quad\underset{t\in[0,1]}{max}Area\big(\rho(t)\big)$.
Considering the energy, similarly define
$\mathcal{W}_{E}=\underset{(\rho, \tau)\in[(\beta,
\tau_{0})]}{inf}\quad\underset{t\in[0,1]}{max}E\big(\rho(t),
\tau(t)\big)$ \footnote{The Teichm$\ddot{u}$ller space
$\mathcal{T}_{1}$ is simply connected, so we do not need to consider
the homotopy class of conformal structures, i.e.$[(\rho, \tau)]$ is
the same as $[\rho]$.}.
\end{Definition}

In fact, we will show that
$\mathcal{W}=\mathcal{W}_{E}$ in Remark \ref{equivalence of minmax
area and energy}. What we are interested is the case when $\mathcal{W}>0$.
So we assume that $\mathcal{W}>0$ in the following.


\subsection{Sketch of the variational approach}\label{sketch of the variational
approach}

\textbf{Question}: Whether one can find a minimal torus or a minimal
torus together with several minimal spheres with total area equal
$\mathcal{W}$? Here we will follow the method of Colding and
Minicozzi. We want to reduce the variational problem for the area
functional to that of the energy functional, i.e. to change a
variational problem in $\Omega$ to one in $\tilde{\Omega}$.
\textbf{Firstly} choose a sequence
$\tilde{\gamma}_{n}(t)\in[\beta]$, such that
$\lim_{n\rightarrow\infty}\underset{t\in[0,1]}{max}Area\big(\tilde{\gamma}_{n}(t)\big)=\mathcal{W}$.
By a smoothing argument, we can assume $\tilde{\gamma}_{n}(t)$
varies in the $C^{2}$ class w.r.t $t$, i.e.
$\tilde{\gamma}_{n}(t)\in C^{0}\Big([0,1],C^{2}(T^{2}_{0}, N)\Big)$.
Pull back the ambient metric
$\tilde{g}_{n}(t)=\tilde{\gamma}_{n}(t)^{*}h$. We want to show that
$\tilde{g}_{n}(t)$, which may be degenerate, determine a family of
marks $\tau_{n}(t)\in\mathcal{T}_{1}$, such that there exist almost
conformal parametrizations $h_{n}(t): T^{2}_{\tau_{n}(t)}\rightarrow
T^{2}_{\tilde{g}_{n}(t)}$ isotopic to $i_{\tau_{n}(t)}$. Hence the
reparametrization $\big(\gamma_{n}(t),
\tau_{n}(t)\big)=\Big(\tilde{\gamma}_{n}\big(h_{n}(t), t\big),
\tau_{n}(t)\Big)\in\big[\big(\tilde{\gamma}_{n}(t), \tau_{0}\big)\big]$ have energy
close to area, i.e. $E\big(\gamma_{n}(t),
\tau_{n}(t)\big)-Area\big(\gamma_{n}(t)\big)\rightarrow 0$.
\textbf{Next} we want to perturb $\gamma_{n}(t)$ to $\rho_{n}(t)$ to
get bubble compactness. Clearly, we can not globally change the
mappings on each path to harmonic or almost harmonic ones like in
the Plateau Problem. Local harmonic replacement is a good choice
here, and this is just what Colding and Minicozzi did.
\textbf{Finally} we will study what we will get when the the
corresponding marks $\{\tau_{n}\}\subset\mathcal{T}_{1}$ converge or
degenerate. If the marks $\tau_{n}$ being considered will not
degenerate, we will get a good solution to this variational problem.
In fact, we will show that $\big(\rho_{n}(t), \tau_{n}(t)\big)$ are
almost conformal when their energy are closed to the min-max value
$\mathcal{W}_{E}$.\\

We will give details of the above approach in the following
sections.


\section{Conformal parametrization}

We will do almost conformal reparametrization for the minimizing
sequence of paths $\tilde{\gamma}_{n}(t)$, and we can assume that
$\tilde{\gamma}_{n}(t)$ have some regularity.

\begin{Lemma}\label{mollifying}
(Lemma D.1 of \cite{CM})Suppose $\tilde{\gamma}_{n}(t)$ are chosen
as a minimizing sequence of paths as above, we can perturb them to
get a new minimizing sequence in the same homotopy class $[\beta]$.
If denoting them still as $\tilde{\gamma}_{n}(t)$, we have
$\tilde{\gamma}_{n}(t)\in C^{0} \Big([0,1],C^{2}(T^{2}_{0},
N)\Big)$.
\end{Lemma}


\subsection{Uniformization for torus}

We need the following uniformization result. For a marked torus
$T^{2}_{\tau}$, we have a standard covering map
$\pi_{\tau}:\mathbb{C}\rightarrow T^{2}_{\tau}$, which is just the
map quotient by the lattices generated by $\{1, \tau\}$. We denote
$\pi_{0}=\pi_{\sqrt{-1}}$.

\begin{Prop}\label{uniformization for torus}
Let $g$ be a $C^{1}$ metric on $T^{2}_{0}$. We can view $g$ as a
metric on the complex plane $\mathbb{C}$, with double periods. Then
there is a unique mark $\tau\in\mathcal{T}_{1}$, and a unique
orientation preserving $C^{1,\frac{1}{2}}$ conformal diffeomorphism
$h:T^{2}_{\tau}\rightarrow T^{2}_{g}$,  such that $h$ is isotopic to
$i_{\tau}$, with normalization that if pulling the map back to
$\mathbb{C}$ by $\pi_{\tau}$ and $\pi_{0}$, it maps $0$ to $0$, $1$
to $1$ and $\tau$ to $\sqrt{-1}$. Furthermore, if $g(t)$ is a family
of $C^{1}$ metrics on $T^{2}_{0}$ which varies continuously in the
$C^{1}$ class, i.e. $g(t)\in C^{1}\big([0,1],\ C^{1}\ metrics\big)$,
and $g(t)\geq\epsilon g_{0}$ for some uniform $\epsilon>0$, let
$\tau(t), h(t)$ be the corresponding marks and normalized conformal
diffeomorphisms, then $\tau(t)$ varies continuously in
$\mathcal{T}_{1}$ and $h(t)$ varies continuously in $C^{0}\cap
W^{1,2}(T^{2}_{\tau(t)}, T^{2}_{0})$.
\end{Prop}

\begin{Remark}
Here the space $C^{0}\cap W^{1,2}(T^{2}_{\tau(t)}, T^{2}_{0})$ have
different domain spaces $T^{2}_{\tau(t)}$, and the continuity is
defined as the Section 2.
\end{Remark}

\begin{Proof}
The existence of a lattice $\{1, \tau\}$ and the conformal
homeomorphism $h:T^{2}_{\tau}\rightarrow T^{2}_{g}$ follows from
Theorem 3.3.2 of \cite{J1} by variational methods.\\

We firstly give the existence of a conformal homeomorphism
satisfying the above normalization. Let $f: T^{2}_{g}\rightarrow
T^{2}_{\tau}$ be the inverse mapping of the conformal homeomorphism
$h$ given by the variational methods. Pulling back $T^{2}_{g}$ to
$\mathbb{C}$ by $\pi_{0}$, $g$ can be viewed as double periodic
metrics $(g_{ij})$. By Lemma \ref{domain metric}, we can write
$g=\lambda|dz+\mu d\overline{z}|^{2}$, with $|\mu|\leq k<1$. Let
$\tilde{f}$ be the lifting of $f$ to the covering space
$\tilde{f}:\mathbb{C}\rightarrow\mathbb{C}$ by $\pi_{0}$ and
$\pi_{\tau}$. After possibly composing with a conformal
diffeomorphism of $T^{2}$, we can assume $\tilde{f}(1)=1$. By the
uniqueness of $\mu$-conformal homeomorphisms which fix
$(0,1,\infty)$ as described in section 6.1, we know that $\tilde{f}$
is just the map $w^{\mu}$ given by Ahlfors and Bers in \cite{AB}.
Since $\tilde{f}$ is orientation preserving,
$\tilde{f}(\sqrt{-1})\in\mathbb{H}$. Denoting
$\tau^{\prime}=\tilde{f}(\sqrt{-1})$, since $f_{\#}$ is
homeomorphism between $\pi_{1}(T^{2}_{0})$ and
$\pi_{1}(T^{2}_{\tau})$, we know $\{1, \tau^{\prime}\}$ is another
generator of the lattice generalized by $\{1, \tau\}$. After pulling
down $\tilde{f}$ by $\pi_{0}$ and $\pi_{\tau^{\prime}}$, we get
$f^{\prime}$. In fact $f^{\prime}$ differs from $f$ by an
automorphism $\pi_{\tau^{\prime}}\circ\pi_{\tau}^{-1}$ of
$T^{2}_{0}$. $f^{\prime}$ maps $T^{2}_{g}$ conformally and
homeomorphicly to $T^{2}_{\tau^{\prime}}$. Since $\tilde{f}$ maps
$1$ to $1$ and $\sqrt{-1}$ to $\tau^{\prime}$, we know that
$f^{\prime}$ is homotopic to $i_{\tau}^{-1}$ by Lemma 2.7.1 of
\cite{J2}. So $f^{\prime}$ and $\tau^{\prime}$ are our unique
conformal homeomorphism and mark, and we will denote them by $f$ and
$\tau$. Let $h=f^{-1}:T^{2}_{\tau}\rightarrow T^{2}_{g}$ be our
unique conformal homeomorphism, then $h$ is isotopic to
$i_{\tau}$.\\

The uniqueness under the above normalization and the continuous
dependence of the conformal homeomorphisms and the marks on the
variance of the metric follow from Appendix \ref{apeendix1}. For a
family of metrics $g(t)$,
$g(t)=\lambda(z)|dz+\mu(t)d\overline{z}|^{2}$, with $|\mu(t)|\leq
k(\epsilon)<1$. Here $\mu(t)$ are double periodic functions on
$\mathbb{C}$ with periods generalized by $\{1, \tau_{0}\}$, and
$\mu(t)=\mu_{t}$ change continuously in the $C^{1}$ class w.r.t $t$
by Lemma \ref{domain metric} and the following Remark \ref{remark of
domain metric}. Let $f(t)$ be the inverse of $h(t)$, with
$\tilde{f}(t)$ and $\tilde{h}(t)$ being pulled back by $\pi_{0}$ and
$\pi_{\tau(t)}$. Hence $\tilde{f}(t)=w^{\mu_{t}}$ are just the maps
given by Ahlfors and Bers described in Appendix \ref{apeendix1}.\\

We will show that $\tau(t)$ vary continuously w.r.t $t$. We know
that $\tau(t)=w^{\mu_{t}}(\sqrt{-1})$, and then
$w^{\mu_{t}}(\sqrt{-1})\rightarrow w^{\mu_{t_{0}}}(\sqrt{-1})$ as
$t\rightarrow t_{0}$. This is because we have convergence under
sphere distance in Lemma \ref{cont1}, i.e.
$d_{S^{2}}\big(w^{\mu_{t}}(\sqrt{-1}),
w^{\mu_{t_{0}}}(\sqrt{-1})\big)\rightarrow 0$. And we know from the
variational methods that $w^{\mu_{t_{0}}}(\sqrt{-1})=\tau_{t_{0}}$
is away from $\infty$, so all $\tau_{t}=w^{\mu_{t}}(\sqrt{-1})$ are
away from $\infty$. Since the sphere distance is equivalent to plane
distance of $\mathbb{C}$, we know
$|w^{\mu_{t}}(\sqrt{-1})-w^{\mu_{0}}(\sqrt{-1})|\rightarrow 0$, i.e.
$\tau(t)\rightarrow\tau(t_{0})$ in $\mathcal{T}_{1}$.\\

We will give the continuous dependence of $h_{t}=f_{t}^{-1}$ on $t$.
The lifting are $\mu(t)$-conformal
$\tilde{h}_{t}:\mathbb{C}_{dwd\overline{w}}\rightarrow\mathbb{C}_{|dz+\mu(t)d\overline{z}|^{2}}$.
Here, we only need to consider $\tilde{h}_{t}$ as mappings defined
on a large ball $B_{R}$, which contains all the parallelograms of
$\{1, \tau_{t}\}$. This is because $\tau_{t}$ vary continuously, so
they will lie on a large ball $B_{R}$ for all $t\in[0,1]$. Here
$\tilde{h}(t)$ are the conformal homeomorphism solutions of Lemma
\ref{cont2}. We know the convergence under sphere distance, i.e.
equation \ref{sphere convergence for the inverse of u-conformal
map}. The image $\tilde{h}(t)(B_{R})$ are restrained to a
neighborhood of $[0,1]\times[0,1]$, since $\tilde{h}(t)$ have
uniform H$\ddot{o}$lder continuity and map parallelograms $\{1,
\tau_{t}\}$ homeomorphicly to $T^{2}_{0}$. So
$\|\tilde{h}_{t}-\tilde{h}_{t_{0}}\|_{L^{\infty}(B_{R})}\rightarrow
0$, as $t\rightarrow t_{0}$, and hence:
\begin{equation}
\|h_{t}-h_{t_{0}}\|_{C^{0}(T^{2}_{\tau_{t}}, T^{2}_{0})}\rightarrow
0.
\end{equation}

From the second convergence in Lemma \ref{cont2}, we know
$\|(\tilde{h}_{t}-\tilde{h}_{t_{0}})_{w}\|_{L^{p}(B_{R})}\rightarrow
0$, as $t\rightarrow t_{0}$, so
$\|(h_{t}-h_{t_{0}})_{w}\|_{L^{p}(T^{2}_{\tau_{t}},
T^{2}_{0})}\rightarrow 0$, and hence:
\begin{equation}
\|h_{t}-h_{t_{0}}\|_{W^{1,2}(T^{2}_{\tau_{t}},
T^{2}_{0})}\rightarrow 0.
\end{equation}
\end{Proof}


\subsection{Construction of the conformal reparametrization}

As above, we consider $\tilde{g}_{n}(t)=\tilde{\gamma}_{n}(t)^{*}h$,
which vary continuously in the $C^{1}$ class. Since there may be
degenerations, we let $g_{n}(t)=\tilde{g}_{n}(t)+\delta_{n} g_{0}$,
where $g_{0}$ is the standard metric of $T^{2}_{0}$, and
$\delta_{n}$ arbitrarily small. The corresponding marks in
$\mathcal{T}_{1}$ and conformal diffeomorphisms are $\tau_{n}(t)$
and $h_{n}(t)$ given by Proposition \ref{uniformization for torus}.
We have the following result.

\begin{Them}\label{conformal parametrization}
Using the above notion, we have reparametrizations
$\big(\gamma_{n}(t), \tau_{n}(t)\big)\in\tilde{\Omega}$ for
$\tilde{\gamma}_{n}(t)$, i.e.
$\gamma_{n}(t)=\tilde{\gamma}_{n}\big(h_{n}(t), t\big)$, such that
$\gamma_{n}(t)\in\big[\tilde{\gamma}_{n}\big]$. And
\begin{equation}\label{equation of conformal parametrization}
E\big(\gamma_{n}(t),
\tau_{n}(t)\big)-Area\big(\gamma_{n}(t)\big)\rightarrow 0,
\end{equation}

as $\delta_{n}\rightarrow 0$.
\end{Them}

\begin{Proof}
We know that $h_{n}(t):T^{2}_{\tau_{n}(t)}\rightarrow
T^{2}_{g_{n}(t)}$ are conformal diffeomorphisms. Let
$\gamma_{n}(t)=\tilde{\gamma}_{n}\big(h_{n}(t),
t\big):T^{2}_{\tau_{n}(t)}\rightarrow N$ be the composition of our
test path with the almost conformal parametrization, we know
$\gamma_{n}(t)\in\Omega$. The continuity of
$t\rightarrow\gamma_{n}(t)$ from $[0,1]$ to $C^{0}\cap
W^{1,2}(T^{2}_{\tau_{n}(t)}, N)$ follows from the continuity of
$t\rightarrow\tilde{\gamma}_{n}(t)$ in $C^{2}$ by Lemma
\ref{mollifying}, and $t\rightarrow h_{n}(t)$ in $C^{0}\cap W^{1,2}$
by Proposition \ref{uniformization for torus}. We will show that
$\gamma_{n}(t)\in[\tilde{\gamma}_{n}]$. From our discussion of
homotopy equivalence of mappings defined on different domains in
Section 2, we view $\gamma_{n}(t)$ as mappings defined on
$T^{2}_{0}$ by composing with
$i_{\tau_{n}(t)}^{-1}:T^{2}_{0}\rightarrow T^{2}_{\tau_{n}(t)}$ and
compare it to $\tilde{\gamma}_{n}(t)$. Since $h_{n}$ are homotopic
equivalent to $i_{\tau_{n}(t)}$ by Proposition \ref{conformal
parametrization}, $h_{n}(t)\circ i_{\tau_{n}(t)}^{-1}$ is homotopic
equivalent to identity map of $T^{2}_{0}$. While $\gamma_{n}$ are
composition of $\tilde{\gamma}_{n}$ with $h_{n}(t)$,
$\gamma_{n}\circ i_{\tau_{n}}^{-1}$ is homotopic equivalent to
$\tilde{\gamma}_{n}$, hence $\gamma_{n}\sim\tilde{\gamma}_{n}$.\\

We can get estimates as in Appendix D of \cite{CM}:
\begin{equation}
\begin{split}
E\big(\gamma_{n}(t), \tau_{n}(t)\big)
&=E\big(h_{n}(t):T^{2}_{\tau_{n}(t)}\rightarrow
T^{2}_{\tilde{g}_{n}(t)}\big)\leq
E\big(h_{n}(t):T^{2}_{\tau_{n}(t)}\rightarrow
T^{2}_{g_{n}(t)}\big)\\
&=Area\big(h_{n}(t):T^{2}_{\tau_{n}(t)}\rightarrow
T^{2}_{g_{n}(t)}\big)\\
&=Area\big(T^{2}_{g_{n}(t)}\big)=\int_{T^{2}_{0}}[det\big(g_{n}(t)\big)]^{\frac{1}{2}}dvol_{0}\\
&=\int_{T^{2}_{0}}[det\big(\tilde{g}_{n}(t)\big)+\delta_{n}Tr_{g_{0}}\tilde{g}_{n}(t)+C(\tilde{g}_{n}(t))\delta_{n}^{2}]^{\frac{1}{2}}dvol_{0}\\
&\leq Area(T^{2}_{\tilde{g}_{n}(t)})+C(\tilde{g}_{n}(t))\sqrt{\delta_{n}}\\
&=Area\big(\gamma_{n}(t):T^{2}_{0}\rightarrow
N\big)+C(\tilde{\gamma}_{n})\sqrt{\delta_{n}}.
\end{split}
\end{equation}

The first and last equality follow from the definition of energy and
area integral, and the second inequality is due to the fact
$\tilde{g}_{n}(t)\leq g_{n}(t)$. Hence we have equation
\ref{equation of conformal parametrization}, as
$\delta_{n}\rightarrow 0$.
\end{Proof}

\begin{Remark}\label{equivalence of minmax area and energy}
We point out that the above Lemma implies that
$\mathcal{W}=\mathcal{W}_{E}$. Since we always have that
$Area(u)\leq E(u, \tau)$, we get $\mathcal{W}\leq\mathcal{W}_{E}$.
We will be done if we know $\mathcal{W}_{E}\leq\mathcal{W}$. By
definition
$\mathcal{W}_{E}\leq\underset{t\in[0,1]}{max}E\big(\gamma_{n}(t),
\tau_{n}(t)\big)$. Since
$\mathcal{W}=\underset{n\rightarrow\infty}{\lim}\underset{t\in[0,1]}{max}Area\big(\gamma_{n}(t)\big)$,
we have
$\mathcal{W}_{E}\leq\underset{n\rightarrow\infty}{\lim}\underset{t\in[0,1]}{max}Area\big(\gamma_{n}(t)\big)=\mathcal{W}$.
\end{Remark}

Now we have reduced the problem in $\Omega$ to that in
$\tilde{\Omega}$ as we discussed above, and it is now easy to deal with
energy $E$ by analytical methods.


\section{Compactification for mappings}

In this case, we can view $\gamma_{n}(t)$ as double periodic
mappings on $\mathbb{C}$, with periods generated by lattices $\{1,
\tau_{n}(t)\}$. So all the mappings have the same domain, but with
different periods, with periods varying continuously. We can do
similar perturbation procedure as what Colding and Minicozzi did in
the case of sphere in \cite{CM}.

\begin{Lemma}\label{compactification}
Let $[\beta]$ and $\mathcal{W}_{E}$ be as in section 2.
For any $\big(\gamma(t), \tau(t)\big)\in[\beta]\subset\tilde{\Omega}$ with
$\underset{t\in[0,1]}{max}E\big(\gamma(t), \tau(t)\big)-
\mathcal{W}_{E}\ll 1$, if $\big(\gamma(t), \tau(t)\big)$ is not
harmonic unless $\gamma(t)$ is a constant map, we can perturb
$\gamma(t)$ to $\rho(t)$, such that $\rho(t)\in[\gamma ]$ and
$E\big(\rho(t), \tau(t)\big)\leq E\big(\gamma(t), \tau(t)\big)$, and
for any $t$ such that
$E\big(\gamma(t),\tau(t)\big)\geq\frac{1}{2}\mathcal{W}_{E}$,
$\rho(t)$ satisfy:

(*) For any finite collection of disjoint balls
$\underset{i}{\cup}B_{i}$ on $T^{2}_{\tau_{t}}$, which can also be
viewed as disjoint balls on the parallegram generated by $\{1,
\tau(t)\}\subset\mathbb{C}$, such that $E\big(\rho(t),
\underset{i}{\cup}B_{i}\big)\leq\epsilon_{0}$, if we let $v$ be the
energy minimizing harmonic map with the same boundary value as
$\rho(t)$ on $\frac{1}{8}\underset{i}{\cup}B_{i}$, then we have:

\begin{equation}\label{compactification formula}
\int_{\frac{1}{8}\underset{i}{\cup}B_{i}}|\nabla\rho(t)-\nabla
v|^{2}\leq\Psi\Big(E\big(\gamma(t), \tau(t)\big)-E\big(\rho(t),
\tau(t)\big)\Big).
\end{equation}

Here $\epsilon_{0}$ is some small constant, and $\Psi$ is a positive
continuous function with $\Psi(0)=0$.
\end{Lemma}

\begin{Remark}
In the paper \cite{CM} of Colding and Minicozzi, all the results
about harmonic maps on disks are still valid here. The other two
most important ingredients are continuity of local maps and
comparison of energy of local harmonic replacements. For the first
one, since all the balls $\underset{i}{\cup}B_{i}$ can be viewed as
balls on $\mathbb{C}$, and $\gamma(t)$ are continuous as mappings on
$\mathbb{C}$, so continuity of $\gamma(t)$ restricted to local balls
is valid. The comparison results are just for a fixed mapping
$\gamma(t)$, and when $t$ is fixed, all the comparison results can
be viewed as on the plane, so we can show that they are still valid
here.
\end{Remark}

We will give the proof by combining results in the following
sections by following the proof of Theorem 2.1 of \cite{CM}. To do
such compactification, we use repeated local harmonic replacements,
which means that we replace the map $u$ on a ball $B$ by the energy-minimizing
map $H(u)$ with the same boundary value as $u$.


\subsection{Harmonic replacement on disks}

In this section, we will list some results about harmonic
replacement on disks with small energy as given in Section 3 of
\cite{CM}. Firstly we recall that for small energy harmonic map,
energy gap can control the difference of $W^{1,2}$-norm. Here
$B_{1}\in\mathbb{R}^{2}$ is the unit disk, and $N$ is the ambient
manifold.

\begin{Them}\label{energy gap}
(Theorem 3.1 of \cite{CM}) There exists a small constant
$\epsilon_{1}$(depending on $N$) such that for all maps $u, v\in
W^{1,2}(B_{1}, N)$ , if $v$ is weakly harmonic with the same
boundary value as $u$, and $v$ has energy less than $\epsilon_{1}$,
then we have:

\begin{equation}
\int_{B_{1}}|\nabla u|^{2}-\int_{B_{1}}|\nabla
v|^{2}\geq\frac{1}{2}\int_{B_{1}}|\nabla u-\nabla v|^{2}.
\end{equation}
\end{Them}

\begin{Remark}
This theorem tells us that for small energy harmonic map, we can use
the gap of energy to control the difference of $W^{1,2}$ norm. Hence
we will focus on the energy gaps when we do harmonic replacement. It
also implies the uniqueness of small energy weakly harmonic map
among maps with the same boundary values(Corollary 3.3 of
\cite{CM}).
\end{Remark}

Using this theorem and boundary regularity of harmonic maps(i.e.
\cite{Q}), we have the following continuity property of harmonic
replacements.

\begin{Cor}\label{continuity of harmonic replacement}
(Corollary 3.4 of \cite{CM}) Let $\epsilon_{1}$ be as in the
previous theorem. Suppose $u\in C^{0}(\overline{B}_{1})\cap
W^{1,2}(B_{1})$ with energy $E(u)\leq\epsilon_{1}$, then there
exists a unique energy minimizing harmonic map $v\in
C^{0}(\overline{B}_{1})\cap W^{1,2}(B_{1})$ with the same boundary
value as $u$. Set $\mathcal{M}=\{u\in C^{0}(\overline{B}_{1})\cap
W^{1,2}(B_{1}), E(u)\leq\epsilon_{1}\}$. $\exists C$ (depending on
$N$), $\forall u_{1}, u_{2}\in\mathcal{M}$, let $w_{1}, w_{2}$ be
the corresponding energy minimizing maps, and let
$E=E(u_{1})+E(u_{2})$, then we have:
\begin{equation}
|E(w_{1})-E(w_{2})|\leq
C\|u_{1}-u_{2}\|_{C^{0}(\overline{B}_{1})}E+C\|\nabla u_{1}-\nabla
u_{2}\|_{L^{2}(B_{1})}E^{\frac{1}{2}}.
\end{equation}

If we denote $v$ by $H(u)$, the mapping $H:\mathcal{M}\rightarrow
\mathcal{M}$ is continuous w.r.t the norm on
$C^{0}(\overline{B}_{1})\cap W^{1,2}(B_{1})$. Here the norm is the
sum of $C^{0}(\overline{B}_{1})$-norm and $W^{1,2}(B_{1})$-norm.
\end{Cor}

We will need the following extension of the above result:

\begin{Cor}\label{continuity of harmonic replacement2}
Suppose $u_{i}, u$ are defined on a ball $B_{1+\epsilon}$ with
energy less than $\epsilon_{1}$. Suppose $u_{i}\rightarrow u$ in
$C^{0}(\overline{B}_{1+\epsilon})\cap W^{1,2}(B_{1+\epsilon})$.
Choose a sequence $r_{i}\rightarrow 1$, and let $w_{i}, w$ be the
mappings which coincide with $u_{i}, u$ outside $r_{i}B_{1}$ and
$B_{1}$ and are energy minimizing inside $r_{i}B_{1}$ and $B_{1}$
respectively. We have $w_{i}\rightarrow w$ in
$C^{0}(\overline{B}_{1+\epsilon})\cap W^{1,2}(B_{1+\epsilon})$.
\end{Cor}

\begin{Proof}
Firstly we show the following claim:

\textbf{Claim:} Let $\tilde{w}_{i}$ be the energy minimizing map
with the same boundary value as $u$ on $r_{i}B_{1}$, then we have:
$\tilde{w}_{i}\rightarrow w$ in
$C^{0}(\overline{B}_{1+\epsilon})\cap W^{1,2}(B_{1+\epsilon})$.

Since $E(u, B_{1+\epsilon})\leq\epsilon_{1}<\epsilon_{SU}$, with
$\epsilon_{SU}$ the constant given in \cite{SU1}, we know that
$\tilde{w}_{i}$ have uniform inner $C^{2,\alpha}$ bounds on $B_{1}$,
so $\forall r<1$, $\tilde{w}_{i}\rightarrow w^{\prime}$ in
$C^{2,\alpha}(B_{r})$, and $w^{\prime}$ is a harmonic map on
$B_{1}$. By scaling argument, we can show that there are no energy
concentration near the boundary of $B_{1}$. So
$\tilde{w}_{i}\rightarrow w^{\prime}$ in $W^{1,2}(B_{1+\epsilon})$.
We also know from \cite{Q}, as indicated by the proof of Corollary
3.4 of \cite{CM} that $\tilde{w}_{i}$ are equi-continuous near
$\partial(r_{i}B_{1})$ and hence equi-continuous near $\partial
B_{1}$ since $r_{i}\rightarrow 1$. So $\tilde{w}_{i}\rightarrow
w^{\prime}$ in $C^{0}(\overline{B}_{1+\epsilon})$. By the uniqueness
of small energy harmonic map of Corollary 3.3 of \cite{CM}, we know
$w^{\prime}=w$. So the claim holds.\\

Let $v_{i}=\Pi(\tilde{w}_{i}+u_{i}-u)$ which have the same boundary
value as $u_{i}$ and $w_{i}$ on $\partial(r_{i}B_{1})$. Here
$\Pi:N_{\delta}\rightarrow N$ is the nearest point projection
defined on a tubular neighborhood $N_{\delta}$. When $\delta$ is
small enough, we have $|d\Pi|\leq2$. So
$\|v_{i}-\tilde{w}_{i}\|_{W^{1,2}(B_{1+\epsilon})}\rightarrow 0$,
hence $\|v_{i}-w\|_{W^{1,2}(B_{1+\epsilon})}\rightarrow 0$ by our
Claim. By Corollary \ref{continuity of harmonic replacement},
$|E(w_{i})-E(\tilde{w}_{i})|\rightarrow 0$, hence
$|E(w_{i})-E(v_{i})|\rightarrow 0$. By Theorem \ref{energy gap},
$\|w_{i}-v_{i}\|_{W^{1,2}(r_{i}B_{1})}\rightarrow 0$. So:
\begin{equation}
\int_{B_{1+\epsilon}}|\nabla w_{i}-\nabla
w|^{2}=\int_{r_{i}B_{1}}|\nabla w_{i}-\nabla
w|^{2}+\int_{B_{1+\epsilon}\backslash r_{i}B_{1}}|\nabla
u_{i}-\nabla w|^{2}\rightarrow 0.
\end{equation}

The convergence to $0$ of the second part of the last term in the
above is due to $u_{i}\rightarrow u$ and $w=u$ outside $B_{1}$.
Hence $w_{i}\rightarrow w$ in $W^{1,2}(B_{1+\epsilon})$.\\

To show the $C^{0}(B_{1+\epsilon})$ convergence, we know from
similar argument as in the proof of the claim, that $w_{i}$ are
equi-continuous near $\partial B_{1}$ by the equi-continuity of
$u_{i}$. Recall that \cite{SU1} gives uniform inner $C^{2,\alpha}$
for $w_{i}$ on $B_{1}$. We have that every subsequence of $w_{i}$
must have $w_{i}\rightarrow w$ in $C^{0}(\overline{B}_{1+\epsilon})$
possibly after taking a further subsequence. So we get
$C^{0}(\overline{B}_{1+\epsilon})$ continuity.
\end{Proof}

\begin{Remark}\label{remark of continuity of harmonic replacement2}
Since we always work on path of mappings, and we will do harmonic
replacement on balls with continuously varying radii, this result
tells us that harmonic replacements will give us another continuous
path if we do harmonic replacement continuously on the initial path.
We can continuously shrink the radii of the disks on which we do
harmonic replacement to $0$, so the new path given by harmonic
replacement can be continuously deformed to the original one, i.e.
they lie in the same homotopy class.
\end{Remark}


\subsection{A comparison result for repeated harmonic replacement}

In this section, we will extend the comparison result of local
harmonic replacements given in Lemma 3.11 of \cite{CM} to the case
of torus. We will use $\mathcal{B}$ to denote a finite collection of
disjoint balls on the complex plane $\mathbb{C}$. If $\mu\in[0,1]$,
we denote $\mu\mathcal{B}$ by a finite collection of balls with the
same centers as $\mathcal{B}$, but the radii $\mu$ timing those of
$\mathcal{B}$. If $u$ is a $C^{0}\cap W^{1,2}$ mapping on the
complex plane with small energy on a collection $\mathcal{B}$, let
$H(u, \mathcal{B})$ be the mapping which coincides with $u$ outside
$\mathcal{B}$, and is the energy minimizing inside $\mathcal{B}$. If
$\mathcal{B}_{1}, \mathcal{B}_{2}$ are two such collections, we
denote $H(u, \mathcal{B}_{1}, \mathcal{B}_{2})$ to be $H\big(H(u,
\mathcal{B}_{1}), \mathcal{B}_{2}\big)$. We will give the
relationship between the energy gaps of $u$, $H(u, \mathcal{B}_{1})$
and $H(u, \mathcal{B}_{1}, \mathcal{B}_{2})$.

\begin{Lemma}\label{comparison}
Fix a torus $T^{2}_{\tau}$ with mark $\tau\in\mathcal{T}_{1}$, and
$u\in C^{0}\cap W^{1,2}(T^{2}_{\tau}, N)$. Let $\mathcal{B}_{1}$,
$\mathcal{B}_{2}$ be two finite collection of disjoint balls on
$T^{2}_{\tau}$, which can also be viewed as collections of disjoint
balls on $\mathbb{C}$. If $E(u,
\mathcal{B}_{i})\leq\frac{1}{3}\epsilon_{1}$, with $\epsilon_{1}$ as
in Theorem \ref{energy gap} for $i=1,2$, then there exists a
constant $k$ depending on $N$, such that:
\begin{equation}\label{comparison inequality1}
E(u)-E[H(u, \mathcal{B}_{1}, \mathcal{B}_{2})]\geq
k\bigg(E(u)-E[H(u, \frac{1}{2}\mathcal{B}_{2})]\bigg)^{2},
\end{equation}

and for any $\mu\in[\frac{1}{8}, \frac{1}{2}]$,
\begin{equation}\label{comparison inequality2}
\frac{1}{k}\big(E(u)-E[H(u,
\mathcal{B}_{1})]\big)^{\frac{1}{2}}+E(u)-E[H(u,
2\mu\mathcal{B}_{2})]\geq E[H(u, \mathcal{B}_{1})]-E[H(u,
\mathcal{B}_{1}, \mu\mathcal{B}_{2})].
\end{equation}
\end{Lemma}

\begin{Remark}
We know from the energy minimizing property of small energy harmonic
maps that the following estimates hold:
\begin{equation}
E(u)-E[H(u, \mathcal{B}_{1}, \mathcal{B}_{2})]\geq E(u)-E[H(u,
\frac{1}{2}\mathcal{B}_{1})].
\end{equation}

So the above three inequalities tell us the relationship of energy
improvement between any two successive harmonic replacements.
\end{Remark}

We will give the proof by constructing comparison mappings. We will
use the following Lemma in our construction. Let $B_{R}$ be the ball
of radius $R$ and center $0$ in $\mathbb{C}$, and $N$ the ambient
manifold.

\begin{Lemma}\label{construction from boundary}
(Lemma 3.14 of \cite{CM}) There exists a $\delta$ and a large
constant $C$ depending on $N$, such that for any $f,g\in C^{0}\cap
W^{1,2}(\partial B_{R}, N)$, if $f,g$ are equal at some point on
$\partial B_{R}$, and:
\begin{equation}
R\int_{\partial B_{R}}|f^{\prime}-g^{\prime}|^{2}\leq\delta^{2},
\end{equation}

we can find some $\rho\in(0,\frac{1}{2}R]$, and a mapping $w\in
C^{0}\cap W^{1,2}(B_{R}\backslash B_{R-\rho}, N)$ with
$w|_{B_{R}}=f$, $w|_{B_{R-\rho}}=g$, which satisfies estimates:
\begin{equation}
\int_{B_{R}\backslash B_{R-\rho}}|\nabla w|^{2}\leq
C\big(R\int_{\partial
B_{R}}|f^{\prime}|^{2}+|g^{\prime}|^{2}\big)^{\frac{1}{2}}\big(R\int_{\partial
B_{R}}|f^{\prime}-g^{\prime}|^{2}\big)^{\frac{1}{2}}.
\end{equation}
\end{Lemma}

\begin{Remark}
The condition and result of this Lemma are all scaling invariant, so
we can apply it to balls of any radius $R$.
\end{Remark}

\begin{Proof}
(of Lemma \ref{comparison}) We know that both $u$ and $H(u,
\mathcal{B}_{1})$ have energy less than $\frac{2}{3}\epsilon_{2}$ on
$\mathcal{B}_{1}\cup\mathcal{B}_{2}$, so Theorem \ref{energy gap}
shows that energy gaps can control $W^{1,2}-$norm gaps in this case.
Denote balls in $\mathcal{B}_{1}$ by $B^{1}_{\alpha}$, and balls in
$\mathcal{B}_{2}$ by $B^{2}_{j}$.\\

\textbf{Step 1} (inequality \ref{comparison inequality1}): Since if
the second harmonic replacements are done on balls which are
disjoint with the balls of the first step, the comparison is easy.
So we divide the second class of balls into two disjoint
subcollections
$\mathcal{B}_{2}=\mathcal{B}_{2+}\cup\mathcal{B}_{2-}$, where
$\mathcal{B}_{2+}=\{B^{2}_{j}: \frac{1}{2}B^{2}_{j}\subset
B^{1}_{\alpha}\, or\
\frac{1}{2}B^{2}_{j}\cap\mathcal{B}_{1}=\emptyset\}$ for some
$B^{1}_{\alpha}\in\mathcal{B}_{1}$. We know that:
\begin{equation}
E(u)-E[H(u, \frac{1}{2}\mathcal{B}_{2})]=E(u)-E[H(u,
\frac{1}{2}\mathcal{B}_{2+})]+E(u)-E[H(u,
\frac{1}{2}\mathcal{B}_{2-})].
\end{equation}

We will deal with $\mathcal{B}_{2+}$ and $\mathcal{B}_{2-}$
separately. \\

For $\mathcal{B}_{2+}$, we have:
\begin{equation}
\begin{split}
E(u)-E[H(u,
\frac{1}{2}\mathcal{B}_{2+})]&=\sum_{\{\frac{1}{2}B^{2}_{j}\cap\mathcal{B}_{1}=\emptyset\}}\big(E(u)-E[H(u,
\frac{1}{2}B^{2}_{j})]\big)\\
&+\sum_{\{\frac{1}{2}B^{2}_{j}\subset
B^{1}_{\alpha}\}}\big(E(u)-E[H(u, \frac{1}{2}B^{2}_{j})]\big).
\end{split}
\end{equation}

For balls $\frac{1}{2}B^{2}_{j}\cap\mathcal{B}_{1}=\emptyset$, we
get from the minimizing property of small energy harmonic maps that:
\begin{equation}
\begin{split}
E(u)-E[H(u, \frac{1}{2}B^{2}_{j})] &=E[H(u, \mathcal{B}_{1})]-E[H(u,
\mathcal{B}_{1}, \frac{1}{2}B^{2}_{j})]\\
&\leq E[H(u, \mathcal{B}_{1})]-E[H(u, \mathcal{B}_{1}, B^{2}_{j})].
\end{split}
\end{equation}

So, we have:
\begin{equation}
\begin{split}
\sum_{\{\frac{1}{2}B^{2}_{j}\cap\mathcal{B}_{1}=\emptyset\}}\big(E(u)
&-E[H(u,
\frac{1}{2}B^{2}_{j})]\big)\\
&\leq\sum_{\{\frac{1}{2}B^{2}_{j}\cap\mathcal{B}_{1}=\emptyset\}}E[H(u,
\mathcal{B}_{1})]-E[H(u, \mathcal{B}_{1}, B^{2}_{j})]\\
&\leq E[H(u, \mathcal{B}_{1})]-E[H(u,
\mathcal{B}_{1}, \cup_{\frac{1}{2}B^{2}_{j}\cap\mathcal{B}_{1}=\emptyset}B^{2}_{j})]\\
&\leq E(u)-E[H(u, \mathcal{B}_{1}, \mathcal{B}_{2+})].
\end{split}
\end{equation}

For balls $\frac{1}{2}B^{2}_{j}\subset B^{1}_{\alpha}$, we have
$H(u, \mathcal{B}_{1}, \frac{1}{2}B^{2}_{j})=H(u, \mathcal{B}_{1})$,
so
\begin{equation}
\begin{split}
\int_{B^{2}_{j}}|\nabla H(u, \mathcal{B}_{1}, B^{2}_{j})|^{2} &\leq
\int_{B^{2}_{j}}|\nabla [H(u, \mathcal{B}_{1},
\frac{1}{2}B^{2}_{j})]|^{2}=\int_{B^{2}_{j}}|\nabla H(u,
\mathcal{B}_{1})|^{2}\\
&\leq \int_{B^{2}_{j}}|\nabla H(u,
\frac{1}{2}B^{2}_{j})|^{2}.
\end{split}
\end{equation}

Hence:
\begin{equation}
\int_{B^{2}_{j}}|\nabla u|^{2}-\int_{B^{2}_{j}}|\nabla H(u,
\frac{1}{2}B^{2}_{j})|^{2}\leq \int_{B^{2}_{j}}|\nabla
u|^{2}-\int_{B^{2}_{j}}|\nabla H(u, \mathcal{B}_{1},
B^{2}_{j})|^{2}.
\end{equation}

Summarizing all the results of this case, we have,
\begin{equation}
\begin{split}
\int_{\underset{B^{2}_{j}\subset
B^{1}_{\alpha}}{\cup}B^{2}_{j}}|\nabla u|^{2} &-|\nabla H(u,
\frac{1}{2}B^{2}_{j})|^{2} \leq \int_{\underset{B^{2}_{j}\subset
B^{1}_{\alpha}}{\cup}B^{2}_{j}}|\nabla u|^{2}-|\nabla H(u,
\mathcal{B}_{1}, B^{2}_{j})|^{2}\\
&\leq \int|\nabla u|^{2}-|\nabla
u_{1}|^{2}+\int_{\underset{B^{2}_{j}\subset
B^{1}_{\alpha}}{\cup}B^{2}_{j}}|\nabla u_{1}|^{2}-|\nabla H(u,
\mathcal{B}_{1}, B^{2}_{j})|^{2}
\end{split}
\end{equation}

For the first term, by Theorem \ref{energy gap}, we have
$\int|\nabla u|^{2}-|\nabla u_{1}|^{2}\leq\int|\nabla u-\nabla
u_{1}|^{2}\leq 4\Big(E(u)-E(u_{1})\Big)$. For the second term, we
have $E(u_{1})-E[H(u, \mathcal{B}_{1}, \underset{B^{2}_{j}\subset
B^{1}_{\alpha}}{\cup}B^{2}_{j})]\leq E(u)-E[H(u, \mathcal{B}_{1},
\underset{B^{2}_{j}\subset B^{1}_{\alpha}}{\cup}B^{2}_{j})]$.
Combining them together,
\begin{equation}
E(u)-E[H(u, \frac{1}{2}\mathcal{B}_{2+})]\leq C\big(E(u)-E[H(u,
\mathcal{B}_{1}, \mathcal{B}_{2+})]\big).
\end{equation}\\

For the collection $\mathcal{B}_{2-}$, we should consider balls
separately. Specify a ball $B^{2}_{j}$, such that $B^{2}_{j}\cap
B^{1}_{\alpha}\neq\emptyset$ for some
$B^{1}_{\alpha}\in\mathcal{B}_{1}$. Denote $B^{2}_{j}$ by $B_{R}$,
and $u_{1}=H(u, \mathcal{B}_{1})$. We will compare $E[H(u,
\frac{1}{2}B_{R})]$ with $E[H(u_{1}, B_{R})]$. Using simple measure
theory or the Courant-Leabesgue Lemma(Lemma 3.1.1 of \cite{J1}), we
can find a subset of $[\frac{3}{4}R, R]$ with measure
$\frac{1}{36}R$, such that for any $r$ in this subset, we have:
\begin{equation}
\int_{\partial B_{r}}|\nabla u_{1}-\nabla u|^{2}
\leq\frac{9}{R}\int^{R}_{\frac{3}{4}R}\int_{\partial B_{s}}|\nabla
u_{1}-\nabla u|^{2}\leq\frac{9}{r}\int_{B_{R}}|\nabla u_{1}-\nabla
u|^{2},
\end{equation}
\begin{equation}
\int_{\partial B_{r}}|\nabla u_{1}|^{2}+|\nabla u|^{2}
\leq\frac{9}{R}\int^{R}_{\frac{3}{4}R}\int_{\partial B_{s}}|\nabla
u_{1}|^{2}+|\nabla u|^{2}\leq\frac{9}{r}\int_{B_{R}}|\nabla
u_{1}|^{2}+|\nabla u|^{2}.
\end{equation}

By choosing $\epsilon_{1}$ small enough, we can get $r\int_{\partial
B_{r}}|\nabla u_{1}|^{2}+|\nabla u|^{2}\leq\delta^{2}$ and
$r\int_{\partial B_{r}}|\nabla u_{1}-\nabla u|^{2}\leq\delta^{2}$
with $\delta$ as in the above Lemma \ref{construction from
boundary}. Since $\frac{1}{2}B_{R}\cap B^{1}_{\alpha}\neq\emptyset$,
but $\frac{1}{2}B_{R}\nsubseteq B^{1}_{\alpha}$, $u$ and $u_{1}$
must be equal at some point on $\partial B_{r}$. So from Lemma
\ref{construction from boundary}, we can find a
$\rho\in(0,\frac{1}{2}r]$ and a mapping $w\in C^{0}\cap
W^{1,2}(B_{r}\backslash B_{r-\rho})$ with $w|_{\partial
B_{r}}=u_{1}$, $w|_{\partial B_{r-\rho}}=u$, and:
\begin{equation}\label{construction of w}
\begin{split}
\int_{B_{r}\backslash B_{r-\rho}}|\nabla w|^{2} &\leq
C\big(r\int_{\partial B_{r}}|\nabla u_{1}-\nabla
u|^{2}\big)^{\frac{1}{2}}\big(r\int_{\partial B_{r}}|\nabla
u_{1}|^{2}+|\nabla u|^{2}\big)^{\frac{1}{2}}\\
&\leq C\big(\int_{B_{R}}|\nabla u_{1}-\nabla
u|^{2}\big)^{\frac{1}{2}}\big(\int_{B_{R}}|\nabla u_{1}|^{2}+|\nabla
u|^{2}\big)^{\frac{1}{2}}.
\end{split}
\end{equation}

Define a comparison map $v$ on $B_{R}$ such that:
\begin{displaymath}
v = \left\{ \begin{array}{ll}
u_{1} & \textrm{on $B_{R}\backslash B_{r}$}\\
w & \textrm{on $B_{r}\backslash B_{r-\rho}$}\\
H(u, B_{r})(\frac{r}{r-\rho} x) & \textrm{on $B_{r-\rho}$}
\end{array} \right..
\end{displaymath}

We know $E[H(u_{1}, B_{R})]\leq E(v)$ since $H(u_{1}, B_{R})$ is
energy minimizing among all maps with the same boundary value on
$B_{R}$. So we have:
\begin{equation}
\begin{split}
\int_{B_{R}}|\nabla H(u_{1}, B_{R})|^{2} &\leq\int_{B_{R}}|\nabla v|^{2}\\
&=\int_{B_{R}\backslash B_{r}}|\nabla
u_{1}|^{2}+\int_{B_{r}\backslash B_{r-\rho}}|\nabla
w|^{2}+\int_{B_{r-\rho}}|\nabla H(u, B_{r})(\frac{r}{r-\rho}\
\cdot)|^{2}\\
&=\int_{B_{R}\backslash B_{r}}|\nabla
u_{1}|^{2}+\int_{B_{r}\backslash B_{r-\rho}}|\nabla
w|^{2}+\int_{B_{r}}|\nabla H(u, B_{r})|^{2}.
\end{split}
\end{equation}

The second equation is due to conformal invariance of the Dirichlet
integral. Hence
\begin{equation}
\begin{split}
&\int_{\frac{1}{2}B_{R}}|\nabla
u|^{2}-\int_{\frac{1}{2}B_{R}}|\nabla H(u, \frac{1}{2}B_{R})|^{2}
\leq\int_{B_{r}}|\nabla
u|^{2}-\int_{B_{r}}|\nabla H(u, B_{r})|^{2}\\
&\leq\int_{B_{r}}|\nabla u|^{2}-\int_{B_{R}}|\nabla H(u_{1},
B_{R})|^{2}+\int_{B_{r}\backslash B_{r-\rho}}|\nabla
w|^{2}+\int_{B_{R}\backslash B_{r}}|\nabla u_{1}|^{2}\\
&\leq\int_{B_{R}}|\nabla u_{1}|^{2}-\int_{B_{R}}|\nabla H(u_{1},
B_{R})|^{2}+\int_{B_{r}\backslash B_{r-\rho}}|\nabla
w|^{2}\\
&+\int_{B_{r}}|\nabla u|^{2}-\int_{B_{r}}|\nabla u_{1}|^{2}.
\end{split}
\end{equation}

By argument similar to the above, we know $\int|\nabla
u|^{2}-|\nabla u_{1}|^{2}\leq 4\Big(E(u)-E(u_{1})\Big)$. Put the
estimates \ref{construction of w} into the above inequality, and sum
over $B^{2}_{j}\in\mathcal{B}_{2-}$:
\begin{equation}
\begin{split}
E(u)-E[H(u, \frac{1}{2}\mathcal{B}_{2-})]&\leq E(u_{1})-E[H(u_{1},
\mathcal{B}_{2-})]\\
&+C\big(E(u)-E(u_{1})\big)^{\frac{1}{2}}+E(u)-E(u_{1})\\
&=E(u)-E[H(u_{1},
\mathcal{B}_{2-})]+C\big(E(u)-E(u_{1})\big)^{\frac{1}{2}}\\
&\leq E(u)-E[H(u, \mathcal{B}_{1},
\mathcal{B}_{2})]+C\big(E(u)-E[H(u, \mathcal{B}_{1},
\mathcal{B}_{2})]\big)^{\frac{1}{2}}.
\end{split}
\end{equation}

Using the fact that all the maps have energy less than
$\frac{1}{2}\epsilon_{1}$, we have:
\begin{equation}
E(u)-E[H(u, \frac{1}{2}\mathcal{B}_{2-})]\leq
C^{\prime}\big(E(u)-E[H(u, \mathcal{B}_{1},
\mathcal{B}_{2})]\big)^{\frac{1}{2}}.
\end{equation}

Combining results on $\mathcal{B}_{2+}$ and $\mathcal{B}_{2-}$, we
have:
\begin{equation}
E(u)-E[H(u, \frac{1}{2}\mathcal{B}_{2})]\leq C\big(E(u)-E[H(u,
\mathcal{B}_{1}, \mathcal{B}_{2})]\big)^{\frac{1}{2}},
\end{equation}
i.e. the first inequality \ref{comparison inequality1}.\\

\textbf{Step 2} (inequality \ref{comparison inequality2}): In this
step, we also divide $\mathcal{B}_{2}$ into two classes with
$\mathcal{B}_{2+}=\{B^{2}_{j}:\mu B^{2}_{j}\subset B^{1}_{\alpha}\
or\ \mu B^{2}_{j}\cap\mathcal{B}_{1}=\emptyset\}$. For $\mu
B^{2}_{j}\subset B^{1}_{\alpha}$, we have $H(u,
\mathcal{B}_{1})=H(u, \mathcal{B}_{1}, \mu B^{2}_{j})$, so we need
not to consider such ball. For $\mu
B^{2}_{j}\cap\mathcal{B}_{1}=\emptyset$, we have:
\begin{equation}
E[H(u, \mathcal{B}_{1})]-E[H(u, \mathcal{B}_{1}, \mu
B^{2}_{j})]=E(u)-E[H(u, \mu B^{2}_{j})]\leq E(u)-E[H(u, 2\mu
B^{2}_{j})].
\end{equation}

So summing all the balls in $\mathcal{B}_{2+}$, we have:
\begin{equation}
E[H(u, \mathcal{B}_{1})]-E[H(u, \mathcal{B}_{1},
\mu\mathcal{B}_{2+})]\leq E(u)-E[H(u, 2\mu\mathcal{B}_{2+})].
\end{equation}\\

For the class $\mathcal{B}_{2-}$, we use similar method as above.
The difference are that $B_{R}=2\mu B^{2}_{j}$, and in the
definition of $v$, the role of $u$, $u_{1}$ changed:
\begin{displaymath}
v = \left\{ \begin{array}{ll}
u & \textrm{on $B_{R}\backslash B_{r}$}\\
w & \textrm{on $B_{r}\backslash B_{r-\rho}$}\\
H(u_{1}, B_{r})(\frac{r}{r-\rho}\ x) & \textrm{on $B_{r-\rho}$}
\end{array} \right..
\end{displaymath}

So we have:
\begin{equation}
\int_{B_{R}}|\nabla H(u, B_{R})|^{2}\leq\int_{B_{R}\backslash
B_{r}}|\nabla u|^{2}+\int_{B_{r}\backslash B_{r-\rho}}|\nabla
w|^{2}+\int_{B_{r}}|\nabla H(u_{1}, B_{r})|^{2}.
\end{equation}

And
\begin{equation}
\begin{split}
&\int_{\frac{1}{2}B_{R}}|\nabla
u_{1}|^{2}-\int_{\frac{1}{2}B_{R}}|\nabla H(u_{1},
\frac{1}{2}B_{R})|^{2}\leq\int_{B_{r}}|\nabla
u_{1}|^{2}-\int_{B_{r}}|\nabla H(u_{1}, \frac{1}{2}B_{R})|^{2}\\
&\leq\int_{B_{r}}|\nabla u_{1}|^{2}-\int_{B_{R}}|\nabla H(u,
B_{R})|^{2}+\int_{B_{R}\backslash B_{r}}|\nabla
u|^{2}+\int_{B_{r}\backslash B_{r-\rho}}|\nabla w|^{2}\\
&\leq\int_{B_{R}}|\nabla u|^{2}-\int_{B_{R}}|\nabla H(u,
B_{R})|^{2}+\int_{B_{r}\backslash B_{r-\rho}}|\nabla
w|^{2}+\int_{B_{r}}|\nabla u_{1}|^{2}-|\nabla u|^{2}.
\end{split}
\end{equation}

Here we use our argument $\int|\nabla u|^{2}-|\nabla u_{1}|^{2}\leq
4\Big(E(u)-E(u_{1})\Big)$ again. Use estimates \ref{construction of
w} again observing that $u$, $u_{1}$ have local energy less than
$\frac{1}{3}\epsilon_{1}$, and sum over
$B^{2}_{j}\in\mathcal{B}_{2-}$:
\begin{equation}
E(u_{1})-E[H(u_{1}, \mu\mathcal{B}_{2-})]\leq E(u)-E[H(u,
2\mu\mathcal{B}_{2-})]+C\big(E(u)-E(u_{1})\big)^{\frac{1}{2}}.
\end{equation}

Combining results on $\mathcal{B}_{2+}$ and $\mathcal{B}_{2-}$, we
will get inequality \ref{comparison inequality2}.
\end{Proof}


\subsection{Construction of the perturbation}

To construct a perturbation satisfying condition $(*)$ in Lemma
\ref{compactification}, we can reduce to control the energy gaps
instead of $W^{1,2}$-norm. Since we only focus on balls with small
energy, there must be a maximal possible energy decrease for a fixed
map on certain such balls. If we firstly do harmonic replacement on
such balls, we can then control the energy decrease for harmonic
replacement on other small energy balls by the comparison Lemma
\ref{comparison}. For a path $\big(\sigma(t),
\tau(t)\big)\in\tilde{\Omega}$, $\epsilon\in(0, \epsilon_{1}]$,
define: $e_{\epsilon,
\sigma(t)}=\sup_{\mathcal{B}}\{E\big(\sigma(t),
\tau(t)\big)-E[H(\sigma(t), \frac{1}{2}\mathcal{B}), \tau(t)]\}$.
Here $\mathcal{B}$ are chosen as any finite collection of disjoint
balls on $T^{2}_{\tau_{t}}$, satisfying: $E\big(\sigma(t),
\mathcal{B}\big)\leq\epsilon$. We know $e_{\epsilon, \sigma(t)}>0$
if $\big(\sigma(t), \tau(t)\big)$ is not harmonic. $e_{\epsilon,
\sigma}$ has some continuity as follows:

\begin{Lemma}\label{continuity of maximal energy decrease}
Use notations as above, $\forall t\in(0,1)$, if $\sigma(t)$ is not
harmonic, we can find a neighborhood $I^{t}\subset(0,1)$ of $t$
depending on $t$, $\epsilon$ and the path $\sigma$, such that
\begin{equation}
e_{\frac{1}{2}\epsilon, \sigma(s)}\leq 2 e_{\epsilon, \sigma(t)},
\end{equation}
for $s\in 2 I^{t}$.
\end{Lemma}

\begin{Proof}
$\sigma(t)\in C^{0}\cap W^{1,2}(T^{2}_{\tau_{t}})$ can be viewed as
defined on a uniform domain $B_{R}\subset\mathbb{C}$ with $\{1,
\tau(t)\}\subset B_{R}$ for all $t\in[0,1]$, i.e. $\sigma\in
C^{0}\big([0,1], C^{0}\cap W^{1,2}(B_{R}, N)\big)$. Since
$e_{\epsilon, \sigma(t)}>0$, we can find a neighborhood $\tilde{I}$
of $t$ such that for all $s\in\tilde{I}$, and for any
$\mathcal{B}\subset B_{R}$, we have
\begin{equation}
\frac{1}{2}\int_{\mathcal{B}}|\nabla\sigma(s)-\nabla\sigma(t)|^{2}\leq
min\{\frac{1}{4}e_{\epsilon, \sigma(t)}, \frac{1}{2}\epsilon\}.
\end{equation}

For fixed $s\in\tilde{I}$, we can find a finite collection of balls
$\mathcal{B}\subset B_{R}$, such that $E\big(\sigma(s),
\mathcal{B}\big)\leq\frac{1}{2}\epsilon$ and
$E\big(\sigma(s)\big)-E[H(\sigma(s),
\frac{1}{2}\mathcal{B})]\geq\frac{3}{4}e_{\frac{1}{2}\epsilon,
\sigma(s)}$ by the definition of $e_{\frac{1}{2}\epsilon,
\sigma(s)}$. Hence $E\big(\sigma(t), \mathcal{B}\big)\leq
E\big(\sigma(s), \mathcal{B}\big)+\frac{1}{2}\epsilon\leq\epsilon$,
so we have $E\big(\sigma(t)\big)-E[H(\sigma(t),
\frac{1}{2}\mathcal{B})]\leq e_{\epsilon, \sigma(t)}$. Thus:
\begin{equation}
\begin{split}
E\big(\sigma(s)\big)-E[H(\sigma(s), \frac{1}{2}\mathcal{B})] &\leq
|E\big(\sigma(t)\big)-E\big(\sigma(s)\big)|+E\big(\sigma(t)\big)-E[H(\sigma(t),
\frac{1}{2}\mathcal{B})]\\
&+|E[H(\sigma(t), \frac{1}{2}\mathcal{B})]-E[H(\sigma(s),
\frac{1}{2}\mathcal{B})]|.
\end{split}
\end{equation}

By Corollary \ref{continuity of harmonic replacement}, after
possibly shrinking the neighborhood $\tilde{I}$ to a smaller one
$I$, we will have
$|E\big(\sigma(t)\big)-E\big(\sigma(s)\big)|\leq\frac{1}{4}e_{\epsilon,
\sigma(t)}$ and $|E[H(\sigma(t),
\frac{1}{2}\mathcal{B})]-E[H(\sigma(s),
\frac{1}{2}\mathcal{B})]|\leq\frac{1}{4}e_{\epsilon, \sigma(t)}$. So
we know $E\big(\sigma(s)\big)-E[H(\sigma(s),
\frac{1}{2}\mathcal{B})]\leq\frac{3}{2}e_{\epsilon, \sigma(t)}$, and
hence $e_{\frac{1}{2}\epsilon, \sigma(s)}\leq 2e_{\epsilon,
\sigma(t)}$.
\end{Proof}

Now we will find a good family of coverings of the time parameter on
which we do harmonic replacement for fixed $\gamma(t)$. In fact,
there will be at most two overlaps for these coverings for a fixed
time $t$.

\begin{Lemma}
Let $\big(\gamma(t), \tau(t)\big)$ be as in Lemma
\ref{compactification}, $B_{R}\supset\{1, \tau(t)\}$ as above. There
exist $m$ collection of disjoint balls $\mathcal{B}_{1},\cdots,
\mathcal{B}_{m}\subset B_{R}$, which are disjoint balls on
$T^{2}_{\tau(t)}$ after quotient by $\{1, \tau(t)\}$, and continuous
functions $r_{j}:[0,1]\rightarrow[0,1]$, $j=1,\cdots,m$, satisfying:

$1^{\circ}$. At most two $r_{j}$ are positive for a fixed $t$, and
$E\big(\gamma(t),
r_{j}(t)\mathcal{B}_{j}\big)\leq\frac{1}{3}\epsilon_{1}$;

$2^{\circ}$. If $t\in[0,1]$, such that $E\big(\gamma(t),
\tau(t)\big)\geq\frac{1}{2}\mathcal{W}$, there exists a $j$, such
that $E\big(\gamma(t)\big)-E[H(\gamma(t),
\frac{1}{2}r_{j}\mathcal{B}_{j})]\geq\frac{1}{8}e_{\frac{1}{8}\epsilon_{1},
\gamma(t)}$.
\end{Lemma}

\begin{Proof}
By continuity, $I=\{t\in[0,1]:E\big(\gamma(t),
\tau(t)\big)\geq\frac{1}{2}W\}$ is a compact subset of $(0,1)$,
since the boundary maps $\gamma(0), \gamma(1)$ have energy almost
$0$ by our almost conformal parametrization. Since $\gamma(t)$ has
no nonconstant harmonic slices, $\forall t\in I$, we can find a
finite collection of disjoint balls $\mathcal{B}_{t}$, such that,
$E\big(\gamma(t), \mathcal{B}_{t}\big)\leq\frac{1}{4}\epsilon_{1}$,
and:
\begin{equation}
E\big(\gamma(t)\big)-E[H(\gamma(t),
\frac{1}{2}\mathcal{B}_{t})]\geq\frac{1}{2}e_{\frac{1}{4}\epsilon_{1},
\gamma(t)}>0.
\end{equation}

By Lemma \ref{continuity of maximal energy decrease} and continuity
of $\gamma$, we can find a neighborhood $I^{t}\ni t$, such that:
$e_{\frac{1}{8}\epsilon_{1}, \gamma(s)}\leq 2
e_{\frac{1}{4}\epsilon_{1}, \gamma(t)}$, and $E\big(\gamma(s),
\mathcal{B}_{t}\big)\leq\frac{1}{3}\epsilon_{1}$ for $s\in 2I^{t}$.
By the continuity of harmonic replacement Corollary \ref{continuity
of harmonic replacement}, after possibly shrinking $I^{t}$, we can
get for $s\in 2I^{t}$:
\begin{equation}
|\{E\big(\gamma(t)\big)-E[H(\gamma(t),
\frac{1}{2}\mathcal{B}_{t})]\}-\{E\big(\gamma(s)\big)-E[H(\gamma(s),
\frac{1}{2}\mathcal{B}_{t})]\}|\leq\frac{1}{4}e_{\frac{1}{4}\epsilon_{1},
\gamma(t)}.
\end{equation}

So we have $E\big(\gamma(s)\big)-E[H(\gamma(s),
\frac{1}{2}\mathcal{B}_{t})]\geq\frac{1}{4}e_{\frac{1}{4}\epsilon_{1},
\gamma(t)}\geq\frac{1}{8}e_{\frac{1}{8}\epsilon_{1}, \gamma(s)}$,
for $s\in 2I^{t}$. By the compactness of $I$, we can find a finite
covering $\{I^{t_{i}}\}$ of $I$, and we can shrink $I^{t_{i}}$ such
that each $I^{i}$ intersects at most two $I^{t_{k}}$, and these two
intervals do not intersect with each other. Choose
$\mathcal{B}_{j}=\mathcal{B}_{t_{j}}$, and choose $r_{j}$ which are
equal to $1$ on $I^{t_{j}}$, and $0$ outside $2I^{t_{j}}$. We also
urge that $r_{j}(t)=0$, if $t$ lies in other interval $I^{t_{l}}$
which does not intersect with $I^{t_{j}}$. It is easy to see these
$\mathcal{B}_{j}$ and $r_{j}$ satisfy the Lemma.
\end{Proof}

\begin{Proof}
(of Lemma \ref{compactification}) Choose the covering
$\mathcal{B}_{j}$ and functions $r_{j}$ as the above Lemma. Let
$\gamma^{0}(t)=\gamma(t)$, and $\gamma^{k}(t)=H\big(\gamma^{k-1}(t),
r_{k}(t)\mathcal{B}_{k}\big)$, for $k=1, \cdots, m$. and let
$\rho(t)=\gamma^{m}(t)$. We will show that $\rho\in[\gamma]$. By
Corollary \ref{continuity of harmonic replacement2}, we know
$t\rightarrow\gamma^{k}(t)$ is continuous from $[0,1]$ to $C^{0}\cap
W^{1,2}$, so $\rho\in\Omega$. Since we can continuously shrink
$r_{j}$ to $0$, and again Corollary \ref{continuity of harmonic
replacement2} and the Remark \ref{remark of continuity of harmonic
replacement2} show that we can hence continuously deform $\rho$ to
$\gamma$ in $\Omega$. So $\rho\in[\gamma]$. Clearly we have
$E\big(\rho(t)\big)\leq E\big(\gamma(t)\big)$.\\

Now we show property $(*)$. Property $1^{\circ}$ of the above Lemma
shows that there are at most two steps of harmonic replacements from
$\gamma$ to $\rho$, and for fixed $t$ with
$E\big(\gamma(t)\big)\geq\frac{1}{2}\mathcal{W}$ we denote the
possible middle nontrivial harmonic replacement by $\gamma^{k}(t)$.
For any finite collection of disjoint balls
$\mathcal{B}=\underset{i}{\cup}B_{i}$ with $E\big(\rho(t),
\mathcal{B})\leq\frac{1}{12}\epsilon_{1}$, we can assume that
$\gamma(t), \gamma^{k}(t)$ have energy at least
$\frac{1}{8}\epsilon_{1}$ on $\mathcal{B}$, or we have will a lower
bound of $E\big(\gamma(t)\big)-E\big(\rho(t)\big)$, hence inequality
\ref{compactification formula} holds. By property $2^{\circ}$ of the
above Lemma, the energy decrease from $\gamma(t)$ to $\gamma^{k}(t)$
or from $\gamma^{k}(t)$ to $\rho(t)$ is at least
$\frac{1}{8}e_{\frac{1}{8}\epsilon_{1}, \gamma(t)}$. We have
estimates at worst by Lemma \ref{comparison}:
\begin{equation}
E\big(\gamma(t)\big)-E\big(\rho(t)\big)\geq
k\big(\frac{1}{8}e_{\frac{1}{8}\epsilon_{1}, \gamma(t)}\big)^{2}.
\end{equation}

Now using inequality \ref{comparison inequality2} of Lemma
\ref{comparison} with $\mu=\frac{1}{8},\frac{1}{4}$ twice in the
case that two $r_{j}(t)>0$, we have:
\begin{equation}
\begin{split}
E\big(\rho(t)\big)-E[H(\rho(t), \frac{1}{8}\mathcal{B})] &\leq
E\big(\gamma^{k}(t)\big)-E[H(\gamma^{k}(t),
\frac{1}{4}\mathcal{B})]+\frac{1}{k}[E\big(\gamma^{k}(t)\big)-E\big(\rho(t)\big)]^{\frac{1}{2}}\\
&\leq E\big(\gamma(t)\big)-E[H(\gamma(t),
\frac{1}{2}\mathcal{B})]+\frac{1}{k}[E\big(\gamma(t)\big)-E\big(\gamma^{k}(t)\big)]^{\frac{1}{2}}\\
&+\frac{1}{k}[E\big(\gamma(t)\big)-E\big(\rho(t)\big)]^{\frac{1}{2}}\\
&\leq e_{\frac{1}{8}\epsilon_{1}, \gamma(t)}+C
[E\big(\gamma(t)\big)-E\big(\rho(t)\big)]^{\frac{1}{2}}\\
&\leq C[E\big(\gamma(t)\big)-E\big(\rho(t)\big)]^{\frac{1}{2}}.
\end{split}
\end{equation}

It is easy to get similar estimates in the case only one
$r_{j}(t)>0$. If we choose $\epsilon_{0}=\frac{1}{12}\epsilon_{1}$
and $\Psi$ a square root function, together with Theorem \ref{energy
gap}, we will get property $(*)$.
\end{Proof}

\begin{Remark}\label{nonconstant slice of minimizing sequence}
Before going on, we have to give some restrictions on the area
minimizing sequence $\tilde{\gamma}_{n}(t)$. In fact, we can assume
that $\tilde{\gamma}_{n}(t)$ have no non-constant harmonic slices,
i.e. $\big(\tilde{\gamma}_{n}(t), T^{2}_{0}\big)$ is not harmonic
unless it is a constant map. We can do this by a reparametrization
on $T^{2}_{0}$ as on page 10 of \cite{CM}. In fact, we can assume
$\tilde{\gamma}_{n}(t)$ is a constant map on a small region on
$T^{2}_{0}$ by small perturbation. Since $\gamma_{n}(t)$ differ from
$\tilde{\gamma}_{n}(t)$ by a diffeomorphism from
$T^{2}_{\tau_{n}(t)}$ to $T^{2}_{0}$, $\gamma_{n}(t)$ is also a
constant map on a small region of $T^{2}_{\tau_{n}(t)}$. Hence
$\gamma_{n}(t)$ is not harmonic unless it is a constant map by the
unique continuation of harmonic maps(Corollary 2.6.1 of \cite{J1}).
So we can apply Lemma \ref{compactification} to $\big(\gamma_{n}(t),
\tau_{n}(t)\big)$. Hence there always exist a min-max sequence
$\big(\rho_{n}(t_{n}), \tau_{n}(t_{n})\big)$, such that
$E\big(\rho_{n}(t_{n}), \tau_{n}(t_{n})\big)\rightarrow \mathcal{W}$
satisfying property $(*)$ of Lemma \ref{compactification}, which
will imply bubbling convergence of $\{\rho_{n}(t_{n}),
\tau_{n}(t_{n})\}$. But we have to remember that we do not know the
behavior of $\tau_{n}(t_{n})$, so we will discuss two cases in the
next section.
\end{Remark}


\section{Convergence results}

In the paper \cite{DLL} of Ding, Li and Liu, they discussed bubbling
convergence results of almost harmonic maps from tori with conformal
structures converging or diverging. If the conformal structures
converge, the sequence of almost harmonic maps will bubbling
converge to a minimal torus together with possibly several minimal
spheres. Here convergence of conformal structures will possibly
ensure existence of a nontrivial minimal torus. If the conformal
structures diverge to infinity, the bubbling limits only contain
several minimal spheres, with the body map from torus degenerate. We
will have similar results for our minimizing sequences
$\big(\rho_{n}(t_{n}), \tau_{n}(t_{n})\big)$. In fact, our sequence
are almost conformal.

\begin{Lemma}\label{almost conformal for perturbed sequence}
If $E\big(\rho_{n}(t_{n}),
\tau_{n}(t_{n})\big)\rightarrow\mathcal{W}$, we have
$E\big(\rho_{n}(t_{n}),
\tau_{n}(t_{n})\big)-Area\big(\rho_{n}(t_{n})\big)\rightarrow 0$.
\end{Lemma}

\begin{Remark}
Although after the perturbation is Section 4, $\big(\rho_{n}(t),
\tau_{n}(t)\big)$ may be far from conformal for some $t\in[0,1]$,
this result tells us that it will still be almost conformal for the
mappings with energy closed to $\mathcal{W}$.
\end{Remark}

\begin{Proof}
We know $\underset{t}{\max} E\big(\gamma_{n}(t),
\tau_{n}(t)\big)\rightarrow\mathcal{W}$, and $E\big(\gamma_{n}(t),
\tau_{n}(t)\big)\geq E\big(\rho_{n}(t), \tau_{n}(t)\big)$. So we
have $E\big(\gamma_{n}(t_{n}),
\tau_{n}(t_{n})\big)-E\big(\rho_{n}(t_{n}),
\tau_{n}(t_{n})\big)\rightarrow 0$. As we know from the construction
from $\gamma_{n}(t)$ to $\rho_{n}(t)$, $\rho_{n}(t)$ is gotten by at
most twice harmonic replacements from $\gamma_{n}(t)$ on balls where
$\gamma_{n}(t)$ have energy less than $\epsilon_{1}$. We denote the
possible middle harmonic replacement by $\gamma^{k}_{n}(t)$ as in
the proof of Lemma \ref{compactification}. From Theorem \ref{energy
gap}, we know that
$\|\nabla\gamma_{n}(t_{n})-\nabla\gamma^{k}_{n}(t_{n})\|_{L^{2}}\leq
4 [E\big(\gamma_{n}(t_{n}),
\tau_{n}(t_{n})\big)-E\big(\gamma^{k}_{n}(t_{n}),
\tau_{n}(t_{n})\big)]\rightarrow 0$, and
$\|\nabla\gamma^{k}_{n}(t_{n})-\nabla\rho_{n}(t_{n})\|_{L^{2}}\leq 4
[E\big(\gamma^{k}_{n}(t_{n}),
\tau_{n}(t_{n})\big)-E\big(\rho_{n}(t_{n}),
\tau_{n}(t_{n})\big)]\rightarrow 0$. Since all the energy of
$\gamma_{n}(t)$, $\rho_{n}(t)$ are bounded, we know that
$|Area\big(\gamma_{n}(t_{n})\big)-Area\big(\rho_{n}(t_{n})\big)|\leq
|Area\big(\gamma_{n}(t_{n})\big)-Area\big(\gamma^{k}_{n}(t_{n})\big)|+|Area\big(\gamma^{k}_{n}(t_{n})\big)-Area\big(\rho_{n}(t_{n})\big)|
\leq
C\{\|\nabla\gamma_{n}(t_{n})-\nabla\gamma^{k}_{n}(t_{n})\|_{L^{2}}+\|\nabla\gamma^{k}_{n}(t_{n})-\nabla\rho_{n}(t_{n})\|_{L^{2}}\}\rightarrow
0$. As $E\big(\gamma_{n}(t_{n}),
\tau_{n}(t_{n})\big)-Area\big(\gamma_{n}(t_{n})\big)\rightarrow 0$,
we have $E\big(\rho_{n}(t_{n}),
\tau_{n}(t_{n})\big)-Area\big(\rho_{n}(t_{n})\big)\rightarrow 0$.
\end{Proof}

To discuss bubble convergence for $\big(\rho_{n}(t_{n}),
\tau_{n}(t_{n})\big)$, we firstly talk about \emph{the convergence
of the metrics} given by $\tau_{n}(t_{n})\in\mathcal{T}_{1}$. In
fact, two metrics $\tau$ and $\tau^{\prime}$ are conformally
equivalent, if they lie in the same orbit of $PSL(2, \mathbb{Z})$.
Denote $\mathcal{M}_{1}=\{z\in\mathbb{C}, |z|\geq 1, Imz>0,
-\frac{1}{2}<Rez\leq\frac{1}{2}, \and\ if\ |z|=1, Rez\geq 0\}$ to be
the fundamental region of $PSL(2, \mathbb{Z})$, which is also the
moduli space of all conformal structures on $T^{2}$. So every such
metric in $\mathcal{T}_{1}$ is conformally equivalent to an element
in $\mathcal{M}_{1}$ after a $PSL(2, \mathbb{Z})$-action. We say
\emph{a sequence $\{\tau_{n}\}$ converge to
$\tau_{0}\in\mathcal{M}_{1}$} if after being conformally translated
to $\{\tau^{\prime}_{n}\}\subset\mathcal{M}_{1}$ by actions in
$PSL(2, \mathbb{Z})$, $\tau^{\prime}_{n}\rightarrow\tau_{0}$. Since
area and energy are all conformally invariant, we can always
consider bubble convergence after conformally changing the domain
metrics to the moduli space $\mathcal{M}_{1}$.\\

There is \emph{a criterion for convergence of conformal structures
on Riemann surfaces} with genus $g$ given by Mumford, i.e. Lemma
3.3.2 in \cite{J1}, or Section 4 in \cite{DLL}. If the lengths of
the shortest closed geodesics on a family of genus $g$ surfaces have
a positive lower bound, then the conformal structures on these
surfaces will converge after possibly taking a subsequence. In the
case of torus $T^{2}$, this criterion is relatively simple. Denote
$\tau=\tau_{1}+\sqrt{-1}\tau_{2}$ to be the conformal structure on a
marked torus, and we use the second normalization as discussed
above\footnote{$Area(\omega_{1}, \omega_{2})=1$}. So
$T^{2}_{\tau}=\{\frac{1}{\sqrt{\tau_{2}}},
\frac{\tau_{1}}{\sqrt{\tau_{2}}}+\sqrt{\tau_{2}}\}$. That the
conformal structure $\tau$ degenerate means
$\tau_{2}\rightarrow\infty$. The length of the shortest closed
geodesic on $T^{2}_{\tau}$ has the same order as
$\frac{1}{\sqrt{\tau_{2}}}$. So the criterion is obvious.\\

\begin{Them}\label{convergence theorem}
Using the above notations, let $\big(\rho_{n}(t), \tau_{n}(t)\big)$
be what we get in the last section by perturbation of
$\big(\gamma_{n}(t), \tau_{n}(t)\big)$ as in Lemma
\ref{compactification} and Remark \ref{nonconstant slice of
minimizing sequence}, then all subsequences $\rho_{n}(t_{n})$ with
$E\big(\rho_{n}(t_{n}),
\tau_{n}(t_{n})\big)\rightarrow\mathcal{W}_{E}$, satisfy:

(*) For any finite collection of disjoint balls
$\underset{i}{\cup}B_{i}$ on $T^{2}_{\tau_{n}(t_{n})}$ such that
$E\big(\rho_{n}(t_{n}),
\underset{i}{\cup}B_{i}\big)\leq\epsilon_{0}$, let $v$ be the
harmonic replacement of $\rho_{n}(t_{n})$ on
$\frac{1}{8}\underset{i}{\cup}B_{i}$. We have:
\begin{equation}
\int_{\frac{1}{8}\underset{i}{\cup}B_{i}}|\nabla\rho_{n}(t_{n})-\nabla
v|^{2}\rightarrow 0.
\end{equation}

Here $\epsilon_{0}$ is the small constant given in Lemma
\ref{compactification}. We have the following two possible cases for
$\Big\{\rho_{n}(t_{n}), \tau_{n}(t_{n})\Big\}$:

(1). If $\tau_{n}(t_{n})\rightarrow\tau_{\infty}$ in the above
sense, then there exist a conformal harmonic map $u:\big(T^{2},
\tau_{\infty}\big)\rightarrow N$, and harmonic spheres $\{u_{i}\}$,
such that $\rho_{n}(t_{n})$ bubble converge to $\big(u, u_{1},
\ldots, u_{l}\big)$, with:
\begin{equation}\label{energy identity 1}
\underset{n\rightarrow\infty}{\lim}E\big(\rho_{n}(t_{n}),
\tau_{n}(t_{n})\big)=E(u, \tau_{\infty})+\underset{i}{\sum}E(u_{i}).
\end{equation}

(2). If $\tau_{n}(t_{n})$ diverge, then there exist only several
harmonic spheres $\{u_{i}\}$, such that $\rho_{n}(t_{n})$ bubble
converge to $\big(u_{1}, \ldots, u_{l}\big)$, with body map
degenerated, and
\begin{equation}\label{energy identity 2}
\underset{n\rightarrow\infty}{\lim}E\big(\rho_{n}(t_{n}),
\tau_{n}(t_{n})\big)=\underset{i}{\sum}E(u_{i}).
\end{equation}
\end{Them}

\begin{Remark}
We point out here that property $(*)$ is invariant when we do
recaling in the bubble process. Property $(*)$ also holds when we
conformally change the metrics $\tau_{n}(t_{n})$ to
$\mathcal{M}_{1}$. These two invariance properties ensure us to use
property $(*)$ in all our proof. For case (1), we can use the
bubbling convergence given by Sacks and Uhlenbeck in \cite{SU1}.
Since the area and energy of this sequence will converge to the same
value $\mathcal{W}=\mathcal{W}_{E}$, the energy identity \ref{energy
identity 1} holds. The bubbling limits are the solution of this
variational problem. For case (2), the length of the shortest closed
geodesics will converge to $0$. So by the argument given by Ding, Li
and Liu in Section 4 of \cite{DLL}, we can tear the torus to a long
cylinder. After some conformal scaling, we can assume the radii of
the cylinders equal $1$. So the sequence of almost harmonic mappings
on long cylinders will converge to a set of harmonic spheres by an
argument given in an un-published note \cite{D} of Ding. Similar
argument as case (1) ensures the energy identity \ref{energy
identity 2}.
\end{Remark}

We will need the following Proposition when we prove identities
\ref{energy identity 1} and \ref{energy identity 2}. We denote
$\mathcal{C}_{r_{1}, r_{2}}$ as a part of the cylinder
$S^{1}\times\mathbb{R}$ with radial coordinates between $r_{1}$ and
$r_{2}$. Clearly $\mathcal{C}_{r_{1}, r_{2}}$ is conformally
equivalent to the annulus $B_{e^{-r_{2}}}\backslash B_{e^{-r_{1}}}$.
Here we have to recall the concept of almost harmonic maps defined
by \cite{CM}. Let $N$ be the ambient manifold. For $\nu >0$, we call
$u\in W^{1,2}(\mathcal{C}_{r_{1}, r_{2}}, N)$ a \textbf{$\nu$-almost
harmonic map}(Definition B.27 in \cite{CM}) if for any finite
collection of disjoint balls $\mathcal{B}$ in the conformally
equivalent annulus $B_{e^{-r_{2}}}\backslash B_{e^{-r_{1}}}$ of
$\mathcal{C}_{r_{1}, r_{2}}$, there is an energy minimizing map
$v:\cup_{\mathcal{B}}\frac{1}{8}B\rightarrow N$ with the same
boundary value as $u$ such that:
\begin{equation}
\int_{\frac{1}{8}\mathcal{B}}|\nabla u-\nabla
v|^{2}\leq\nu\int_{\mathcal{C}_{r_{1}, r_{2}}}|\nabla u|^{2}.
\end{equation}

\begin{Prop}\label{almost harmonic maps on necks}
(Proposition B.29 of \cite{CM}) $\forall \delta>0$, there exist
small constants $\nu>0$, $\epsilon_{2}>0$ and large constant $l\geq
1$ (depending on $\delta$ and $N$), such that for any integer $m$,
if $u$ is a $\nu$-almost harmonic map as defined above on
$\mathcal{C}_{-(m+3)l, 3l}$ with $E(u)\leq\epsilon_{2}$, then:
\begin{equation}
\int_{\mathcal{C}_{-ml, 0}}|u_{\theta}|^{2}\leq
7\delta\int_{\mathcal{C}_{-(m+3)l, 3l}}|\nabla u|^{2}.
\end{equation}

Here we use $(\theta, t)$ as coordinates on $S^{1}\times\mathbb{R}$,
and $u_{\theta}$ means the differentiation w.r.t $\theta$.
\end{Prop}

\begin{Proof}
(of Theorem \ref{convergence theorem}) \textbf{Case (1):} We denote
$\rho_{n}=\rho_{n}(t_{n})$, and let $\tau_{n}\in\mathcal{M}_{1}$ be
the corresponding conformal structure of $\tau_{n}(t_{n})$. We
divide the bubbling convergence into several steps, and we will then
focus on the neck parts.\\

\underline{Step 1.} Since $\tau_{n}\rightarrow\tau_{\infty}$, we can
identify a point $x\in T^{2}_{\tau_{\infty}}$ as on
$T^{2}_{\tau_{n}}$ by viewing it as on the fundamental regions of
lattices $\{1, \tau_{\infty}\}$ and $\{1, \tau_{n}\}$ of
corresponding conformal structures. So for any $x\in
T^{2}_{\tau_{\infty}}$, for a fixed small constant
$\epsilon_{1}<\epsilon_{0}$, we can consider a sequence of energy
concentration radii $r_{n}(x)$ defined as follows:
\begin{equation}
r_{n}(x)=sup\{r>0, E\big(\rho_{n}, B(x, r)\big)\leq\epsilon_{1}\}.
\end{equation}

Such $r_{n}(x)$ exist and are positive. Now we say $x$ is an energy
concentration point if $\lim_{n\rightarrow\infty}r_{n}(x)\rightarrow
0$. If $x$ is an energy concentration point, we have that:
\begin{equation}
\underset{r>0}{inf}\{\lim_{n\rightarrow\infty}E(\rho_{n}, B(x,
r))\}\geq\epsilon_{1}.
\end{equation}

Since our sequence $\rho_{n}$ have uniform bounded energy
$2\mathcal{W}$, we know the number of the energy concentration
points are bounded by $2\mathcal{W}/\epsilon_{1}$. Denote these
points by $\{x_{1}, \cdots, x_{m}\}$. If $x\in
T^{2}_{\tau_{\infty}}\backslash\{x_{1}, \cdots, x_{m}\}$, we can
find a $r(x)>0$ such that $E\big(\rho_{n}, B(x,
r(x))\big)\leq\epsilon_{1}$ for all $n$. and by condition $(*)$,
there exist $v_{n}$ which are the energy minimizing harmonic maps
defined on $\frac{1}{8}B(x, r(x))$ with the same boundary value as
$\rho_{n}$, such that
$\|\rho_{n}-v_{n}\|_{W^{1,2}\big(\frac{1}{8}B(x,
r(x))\big)}\rightarrow 0$. Since $E\big(v_{n}, \frac{1}{8}B(x,
r(x))\big)\leq\epsilon_{1}<\epsilon_{SU}$, we know from \cite{SU1}
that $v_{n}$ have uniform interior $C^{2, \alpha}$-estimates on
$\frac{1}{8}B(x, r(x))$, and hence converge to a harmonic map $u$ on
$\frac{1}{9}B(x, r(x))$ in $C^{2, \alpha}$ after taking a
subsequence. Hence $\rho_{n}\rightarrow u$ in
$W^{1,2}\big(\frac{1}{9}B(x, r(x))\big)$. So for any compact subset
$K\subset T^{2}_{\tau_{\infty}}\backslash\{x_{1}, \cdots, x_{m}\}$,
we can cover them by finite many balls $\frac{1}{9}B(x, r(x))$, and
hence $\rho_{n}\rightarrow u$ in $W^{1,2}(K)$ after taking a
subsequence. Here $u$ is a harmonic map defined on $K$. After
exhausting $T^{2}_{\tau_{\infty}}\backslash\{x_{1}, \cdots, x_{m}\}$
by a sequence of compact sets $K_{i}$, and a diagonal argument, we
know $u$ is a harmonic map on
$T^{2}_{\tau_{\infty}}\backslash\{x_{1}, \cdots, x_{m}\}$, and by
the Theorem 3.6 of removable singularity in \cite{SU1}, we know $u$
extends to a harmonic map on $T^{2}_{\tau_{\infty}}$.\\

\underline{Step 2.} We now see what happens near the energy
concentration points. Fix an energy concentration point $x_{i}$, and
denote $r_{n, i}=r_{n}(x_{i})$. Find a small $r>0$, such that $E(u,
B(x_{i}, r))\leq\frac{1}{3}\epsilon_{1}$. We rescale $\rho_{n}$ on
$B(x_{i}, r_{n, i})$. Define $u_{n, i}=\rho_{n}(x_{i}+r_{n,
i}(x-x_{i}))$. So $B(x_{i}, r_{n, i})$ are now rescaled to $B_{1}$,
and $B(x_{i}, r)$ to $B(0, r/r_{n, i})$. $u_{n, i}$ can be viewed as
defined on balls $B(0, r/r_{n, i})$ with radii converging to
infinity. Since the domains converge to the whole complex plane
$\mathbb{C}$, which is conformal equivalent to the sphere $S^{2}$
without the south pole, we can think $u_{n, i}$ as defined on any
compact subsets of $S^{2}$ away from the south pole for $n$ large
enough. Since the property $(*)$ is conformal invariant, we can do
the first step to $u_{n, i}$. We can find finitely many energy
concentration points $\{x_{i, 1}, \cdots, x_{i, m_{i}}\}\subset
S^{2}\backslash south\ pole$, such that $u_{n, i}$ converge to a
harmonic map $u_{i}$ defined on $S^{2}\backslash south\ pole$ in the
sense of the above step, and hence $u_{i}$ is a harmonic sphere
defined on $S^{2}$ by the Theorem of removable singularity. From our
definition, we know that $E(u_{n, i},
B_{1})=\epsilon_{1}<\epsilon_{SU}$. So $x_{i, j}\in S^{2}\backslash
B_{1}$\footnote{Here $B_{1}$ is a unit ball centered at the north
pole.}. A key point is that the total energy of $u_{n, i}$ on
$S^{2}\backslash \{south\ pole\cup B_{1}\}$ is decreased by a finite
amount $\epsilon_{1}$ compared to the original $u_{n, i}$, as $u_{n,
i}|_{B_{1}}$ taking the energy. We call such rescaling and
convergence procedure bubbling convergence.\\

\underline{Step 3.} We can repeat the bubbling convergence given in
step 2 for $u_{n, i}$ on balls centered at $x_{i, j}$. We point out
here that there are only finite many such steps, and then the
bubbling convergence stops. Each time, we come from a sequence of
maps $u_{n}$ defined on a small ball $B_{r}$, and we rescale them to
exhaust the whole complex plane. Each time $u_{n}|_{B_{1}}$ take a
finite amount of energy after recaling. So after several steps, the
total energy of $u_{n}$ will be less than
$\epsilon_{1}<\epsilon_{SU}$, and there will be no energy
concentration points. The bubbling convergence stops.\\

\underline{Step 4.} We will discuss energy identity \ref{energy
identity 1} now. We can decompose $T^{2}_{\tau_{n}}$ into the bubble
part $\underset{i}{\cup}B(x_{i}, r)$ and the body part
$T^{2}_{\tau_{n}}\backslash\underset{i}{\cup}B(x_{i}, r)$. So the
total energy has decomposition $E\big(\rho_{n},
T^{2}_{\tau_{n}}\big)=E\big(\rho_{n},
T^{2}_{\tau_{n}}\backslash\underset{i}{\cup}B(x_{i},
r)\big)+\underset{i}{\sum}E\big(\rho_{n}, B(x_{i}, r)\big)$. Now we
can calculate the energy of the first limit map $u_{0}$ as follows:
\begin{equation}
E(u_{0})=\lim_{r\rightarrow 0}E\big(u_{0},
T^{2}_{\tau_{\infty}}\backslash\underset{i}{\cup}B(x_{i}, r)\big)
=\lim_{r\rightarrow 0}\lim_{n\rightarrow\infty}E\big(\rho_{n},
T^{2}_{\tau_{n}}\backslash\underset{i}{\cup}B(x_{i}, r)\big).
\end{equation}

So we only need to show that $\lim_{r\rightarrow
0}\lim_{n\rightarrow\infty}\underset{i}{\sum}E\big(\rho_{n},
B(x_{i}, r)\big)=\underset{i}{\sum}E(u_{i})$. Here $u_{i}$ are the
bubble maps. As in the second step, we know that $u_{i}$ are limits
of $u_{n, i}$ on any compact set of $\mathbb{C}$, so we can
calculate the energy of the first bubble map $u_{i}$ as follows:
\begin{equation}
E(u_{i})=\lim_{R\rightarrow\infty}E\big(u_{i},
B(R)\big)=\lim_{R\rightarrow\infty}\lim_{n\rightarrow\infty}E\big(u_{n,
i}, B(R)\big).
\end{equation}

By the conformal invariance of energy, $E\big(u_{n, i},
B(R)\big)=E\big(\rho_{n}, B(x_{i}, r_{n, i}R)\big)$. So we only need
to show that:
\begin{equation}
\lim_{r\rightarrow 0,
R\rightarrow\infty}\lim_{n\rightarrow\infty}E\big(\rho_{n}, B(x_{i},
r)\backslash B(x_{i}, r_{n, i}R)\big)=0.
\end{equation}

We denote the annulus $A(x_{i}, r, r_{n, i}R)=B(x_{i}, r)\backslash
B(x_{i}, r_{n, i}R)$. Since $A(x_{i}, r, r_{n, i}R)$ is conformally
equivalent to a lang cylinder $\mathcal{C}_{r_{1}, r_{2}}$, with
$r_{1}=-\ln(r_{n, i}R)$, $r_{2}=-\ln(r)$, we call such annuli or
such cylinders \textbf{necks}. So what left is to show that there
will be no energy concentration on necks.\\

\underline{Step 5.} We use Proposition \ref{almost harmonic maps on
necks} to show that necks support no energy in our case. We will use
step 1 as an example, and others follow in the same way. Suppose
there is a lower bound for $E\big(\rho_{n}, \mathcal{C}_{r_{1},
r_{2}}\big)$. Since $\rho_{n}$ will converge to $u_{0}$ on any small
annulus centered at $x_{i}$, and $u_{n, i}$ will converge to $u_{i}$
on any large annulus centered at $0$, for fixed $L>0$ we know that
there can be no energy concentration on $A(x_{i}, re^{-L}, r)$ and
$A(x_{i}, r_{n, i}Re^{-L}, r_{n, i}R)$ for $r\rightarrow 0$ and
$R\rightarrow\infty$. Changing to the cylinder, we know there will
be no energy concentration on a region with fixed length towards
boundary of $\mathcal{C}_{r_{1}. r_{2}}$. Now fix a
$\delta=\frac{1}{140}$, and let $\nu$, $\epsilon_{2}$ and $l$ be as
in Proposition \ref{almost harmonic maps on necks}. We can find a
sub-cylinder $\mathcal{C}_{r^{\prime}_{1}, r^{\prime}_{2}}$ with the
distance between boundaries of them converging to $\infty$, i.e.
$d(\partial\mathcal{C}_{r_{1}. r_{2}},
\partial\mathcal{C}_{r^{\prime}_{1},
r^{\prime}_{2}})\rightarrow\infty$, such that $E(\rho_{n},
\mathcal{C}_{r^{\prime}_{1},
r^{\prime}_{2}})=\frac{1}{2}\epsilon_{2}$. We want to show that
$\rho_{n}$ is $\nu$-almost harmonic on $\mathcal{C}_{r^{\prime}_{1},
r^{\prime}_{2}}$ for $n$ large. In fact for any finite collection of
disjoint balls $\mathcal{B}$ on the annulus, $E(\rho_{n},
\mathcal{B})\leq E(\rho_{n}, \mathcal{C}_{r^{\prime}_{1},
r^{\prime}_{2}})\leq\epsilon_{2}$. We can assume
$\epsilon_{2}\leq\epsilon_{1}$, so $\rho_{n}$ satisfy property
$(*)$, i.e.
$\int_{\frac{1}{8}\mathcal{B}}|\nabla\rho_{n}-v|^{2}\rightarrow 0$,
with $v$ the energy minimizing map. Since $E(\rho_{n},
\mathcal{C}_{r^{\prime}_{1}, r^{\prime}_{2}})$ have uniform lower
bound,
$\int_{\frac{1}{8}\mathcal{B}}|\nabla\rho_{n}-v|^{2}\leq\nu\int_{\mathcal{C}_{r^{\prime}_{1},
r^{\prime}_{2}}}|\nabla\rho_{n}|^{2}$ hold for $n$ large enough. We
can assume we first do the above on a cylinder a little bit larger
than $\mathcal{C}_{r^{\prime}_{1}, r^{\prime}_{2}}$, then by
Proposition \ref{almost harmonic maps on necks}, we have:
\begin{equation}
\int_{\mathcal{C}_{r^{\prime}_{1},
r^{\prime}_{2}}}|(\rho_{n})_{\theta}|^{2}\leq\frac{1}{10}\int_{\mathcal{C}_{r^{\prime}_{1},
r^{\prime}_{2}}}|\nabla\rho_{n}|^{2}.
\end{equation}

Hence we have a lower bound on the gap between energy and area.
\begin{equation}
\begin{split}
E(\rho_{n}, \mathcal{C}_{r^{\prime}_{1},
r^{\prime}_{2}})-Area(\rho_{n}, \mathcal{C}_{r^{\prime}_{1},
r^{\prime}_{2}})&=\frac{1}{2}\int_{\mathcal{C}_{r^{\prime}_{1},
r^{\prime}_{2}}}|(\rho_{n})_{t}|^{2}+|(\rho_{n})_{\theta}|^{2}-2|(\rho_{n})_{t}\times(\rho_{n})_{\theta}|\\
&\geq\frac{1}{8}\int_{\mathcal{C}_{r^{\prime}_{1},
r^{\prime}_{2}}}|(\rho_{n})_{t}|^{2}-|(\rho_{n})_{\theta}|^{2}.
\end{split}
\end{equation}

So $E(\rho_{n}, \mathcal{C}_{r^{\prime}_{1},
r^{\prime}_{2}})-Area(\rho_{n}, \mathcal{C}_{r^{\prime}_{1},
r^{\prime}_{2}})$ have a lower bound by the above estimates. It is a
contradiction to $E(\rho_{n}(t),
\tau_{n}(t))-Area(\rho_{n}(t))\rightarrow 0$ given in Lemma
\ref{almost conformal for perturbed sequence}.\\

\textbf{Case (2).} We use $(t, \theta)$ as parameters on
$T^{2}_{\tau_{n}}$. In fact, we assume $\arg(\tau_{n})=\theta_{n}$,
and let
$z^{\prime}=t+\sqrt{-1}\theta=e^{-\sqrt{-1}(\frac{1}{2}\pi-\theta_{n})}z$
be another conformal parameter system on $T^{2}_{\tau_{n}}$. We
conformally expand the torus such that the length of the circle of
parameter $\theta$ is $1$, and the length of parameter $t$ is
denoted by $2l_{n}$. Then we divide the torus $T^{2}_{\tau_{n}}$
into sections with length $1$ in the parameter $t$, i.e.
$T^{2}_{\tau_{n}}=\underset{i}{\cup}S^{1}\times[t_{i}, t_{i+1}]$.\\

We \textbf{claim} that there exists a large $L>0$, such that for $n$
large, there exist $t_{n, 0}$, such that $E(\rho_{n},
S^{1}\times[t_{n, 0}-L, t_{n, 0}+L])>\epsilon_{2}$. If the claim
fails, $\forall L>0$, we can find a subsequence of
$n\rightarrow\infty$, such that $\forall t_{n, i}$, $E(\rho_{n},
S^{1}\times[t_{n, i}-L, t_{n, i}+L])\leq\epsilon_{2}$. After
possibly extending some $[t_{n, i}-L, t_{n, i}+L]$, we have
$E(\rho_{n}, S^{1}\times[t_{n, i}-L, t_{n, i}+L])=\epsilon_{2}$. So
$\rho_{n}|_{[t_{n, i}-L, t_{n, i}+L]}$ satisfy condition of
Proposition \ref{almost harmonic maps on necks}, and hence is a
contradiction to Lemma \ref{almost conformal for perturbed sequence}
as argued in step 5 of case 1.\\

Now consider $\rho_{n}:S^{1}\times[t_{n, 0}-l_{n}, t_{n,
0}+l_{n}]\rightarrow N$. There may be bubbles near $t_{n, 0}$.
Argument as in case 1 shows that $\rho_{n}$ converge to a harmonic
map $u_{1}$ defined on $S^{1}\times\mathbb{R}$ besides some energy
concentration points. $u_{1}$ is nontrivial since $E(\rho_{n},
S^{1}\times[t_{n, 0}-L, t_{n, 0}+L])>\epsilon_{2}$. As
$S^{1}\times\mathbb{R}$ is conformally equivalent to
$S^{2}\backslash north\ and\ south\ pole$, we can extend $u_{1}$ to
a harmonic map on $S^{2}$. We can rescale $\rho_{n}$ near the energy
concentration points, and the rescaled map will converge as we
discussed in Case 1 to several bubble maps $\{u_{1, 1}, \cdots,
u_{1, l_{1}}\}$. Energy identity during these bubbles will follow as
in the last step of Case 1 on each long cylinder.  Now we calculate
the total energy:
\begin{equation}
\begin{split}
\lim_{l\rightarrow\infty}\lim_{n\rightarrow\infty}E(\rho_{n},
S^{1}\times[t_{n, 0}-l, t_{n,
0}+l])&=\lim_{l\rightarrow\infty}E(u_{1}, S^{1}\times[-l,
l])+\underset{i}{\sum}E(u_{1, i})\\
&=E(u_{1}, S^{2})+\underset{i}{\sum}E(u_{1, i}).
\end{split}
\end{equation}

So if $\lim_{l\rightarrow\infty}\lim_{n\rightarrow\infty}E(\rho_{n},
S^{2}\times[-l_{n}, -l]\cup[l, l_{n}])=0$, there will be no other
bubbles except for $\{u_{1}, u_{1,1}, \cdots, u_{1,l_{1}}\}$, and
$\lim_{n\rightarrow}E(\rho_{n})=E(u_{1})+\underset{i}{\sum}E(u_{1,
i})$. i.e energy identity \ref{energy identity 2} holds. If
$\lim_{l\rightarrow\infty}\lim_{n\rightarrow\infty}E(\rho_{n},
S^{2}\times[-l_{n}, -l]\cup[l, l_{n}])>0$, we can consider maps on
the other part of the rescaled torus , i.e. we can find another base
point denoted by $t_{n, 1}$, such that $|t_{n, 1}-t_{n,
0}|\rightarrow\infty$ and $E(\rho_{n}, S^{1}\times[t_{n, 1}-L, t_{n,
1}+L])>\epsilon_{2}$. Consider $\rho_{n}:S^{1}\times[t_{n, 1}-l_{n},
t_{n, 1}+l_{n}]\rightarrow N$. We can repeat the above step and get
another set of harmonic spheres $\{u_{2}, u_{2,1}, \cdots, u_{2,
l_{2}}\}$. Since each bubble is a harmonic sphere and must take a
finite mount of energy by \cite{SU1}, there are only finitely many
such steps. We will get all these harmonic spheres $u_{i}$ and
energy identity \ref{energy identity 2} by summing over all the
steps.
\end{Proof}

\textbf{What is left?} The aim of this method is to find a min-max
minimal torus, but only when the conformal structures do not
degenerate can we get a nontrivial minimal torus. So we do want to
know under what condition does there exist a subsequence
$\Big\{\rho_{n}(t_{n}), \tau_{n}(t_{n})\Big\}$ satisfying condition
$(1)$ in the above theorem.


\section{Appendix 1--a uniformization result}\label{apeendix1}

In this section, we discuss a general uniformization theorem on the
complex plane. We will focus on the continuous dependence of the
conformal diffeomorphisms on the variance of general metrics. Let
$g$ be a Riemannian metric on the complex plane $\mathbb{C}$.

\begin{Lemma}\label{domain metric}
In the complex coordinates $\{z, \overline{z}\}$, we can write
$g=\lambda(z)|dz+\mu(z)d\overline{z}|^{2}$. Here $\lambda(z)>0$, and
$\mu(z)$ is complex function on the complex plane with $|\mu|<1$. If
$g\geq\epsilon dzd\overline{z}$, there exists a $k=k(\epsilon)<1$,
such that $|\mu|\leq k$.
\end{Lemma}

\begin{Remark}\label{remark of domain metric}
The proof is just simple calculation. Hence we can always identify a
plane non-degenerate metric with $|dz+\mu(z)d\overline{z}|^{2}$
conformally. In fact, $\mu$ is a rational function of the components
$g_{ij}(z)$, so if a family $g(t)$ vary continuously in the $C^{1}$
class, the corresponding $\mu(t)$ also vary continuously in the
$C^{1}$ class.
\end{Remark}


\subsection{Results in \cite{AB}}

Let us discuss what Ahlfors and Bers did in \cite{AB}. They gave the
\textbf{existence} and \textbf{uniqueness} of conformal
diffeomorphism $w^{\mu}:\mathbb{C}_{|dz+\mu
d\overline{z}|^{2}}\rightarrow\mathbb{C}_{dwd\overline{w}}$ fixing
three points $(0, 1, \infty)$ for any $L^{\infty}$ function $\mu$
with $|\mu|\leq k<1$. Such maps must satisfy the following equation:
\begin{equation}\label{cr1}
w^{\mu}_{\overline{z}}=\mu(z) w^{\mu}_{z}.
\end{equation}

Define function space $B_{p}(\mathbb{C})=C^{1-\frac{2}{p}}\cap
W^{1,p}_{loc}(\mathbb{C})$, where $p>2$ depends on the bound $k$ of
$\mu$. They showed that $w^{\mu}$ are uniformly bounded in
$B_{p}(\mathbb{C})$ for a uniform bound $k$, and that $w^{\mu}$ vary
continuously in $B_{p}(\mathbb{C})$ while $\mu$ varying continuously
in $L^{\infty}(\mathbb{C})$. Suppose $\mu, \nu \in
L^{\infty}(\mathbb{C})$, and $|\mu|,|\nu|\leq k$, with $k<1$. Let
$w^{\mu}, w^{\nu}$ be the corresponding conformal homeomorphisms,
then:

\begin{Lemma}\label{cont1}
(Lemma 16, Theorem 7, Lemma 17, Theorem 8 of \cite{AB})
\begin{equation}
d_{S^{2}}\big(w^{\mu}(z_{1}), w^{\mu}(z_{2})\big)\leq c
d_{S^{2}}(z_{1}, z_{2})^{\alpha},
\end{equation}
\begin{equation}
\|w^{\mu}_{z}\|_{L^{p}(B_{R})}\leq c(R),
\end{equation}
\begin{equation}
d_{S^{2}}\big(w^{\mu}(z), w^{\nu}(z)\big)\leq C
\|\mu-\nu\|_{\infty},
\end{equation}
\begin{equation}
\|(w^{\mu}-w^{\nu})_{z}\|_{L^{p}(B_{R})}\leq
C(R)\|\mu-\nu\|_{\infty}.
\end{equation}

Here $d_{S^{2}}$ is the sphere distance, which is equivalent to the
plane distance of $\mathbb{C}$ on compact sets.
$\alpha=1-\frac{2}{p}$. All constants are uniformly bounded
depending on $k<1$.
\end{Lemma}

\begin{Remark}
This Lemma comes from estimates of equation \ref{cr1}. Here we use
sphere distance because what we concern is just local properties.
\end{Remark}


\subsection{Similar results}

If we write our metrics conformally on $\mathbb{C}$, what we concern
in our case are the conformal homeomorphisms
$h^{\mu}:\mathbb{C}_{dwd\overline{w}}\rightarrow\mathbb{C}_{|dz+\mu
d\overline{z}|^{2}}$ fixing three points $(0,1, \infty)$, which are
just the inverse mappings of those of Ahlfors and Bers. We also
concern the continuous dependence of $h^{\mu}$ in $C^{0}\cap
W^{1,2}_{loc}(\mathbb{C}, \mathbb{C})$ on the variance of $\mu$ in
$C^{1}(\mathbb{C})$. In fact:
\begin{equation}
h^{\mu}(w)=(w^{\mu})^{-1}(w),
\end{equation}

and our mappings satisfy:
\begin{equation}\label{cr2}
h^{\mu}_{\overline{w}}=-\mu(h^{\mu}(w))\overline{h^{\mu}_{w}}.
\end{equation}

If $\mu_n$, are a sequence of metric coefficients as above, such
that $\|\mu_{n}-\mu\|_{C^{1}}\rightarrow 0$, and $h^{\mu_{n}}$ as
above, we want to have results similar to the above:

\begin{Lemma}\label{cont2}
\begin{equation}\label{sphere convergence for the inverse of u-conformal map}
d_{S^{2}}\big(h^{\mu_{n}}, h^{\mu}\big)\rightarrow 0,
\end{equation}
\begin{equation}\label{local Lp convergence of u-conformal map}
\|(h^{\mu_{n}}-h^{\mu})_{w}\|_{L^{p}(B_{R})}\rightarrow 0.
\end{equation}
\end{Lemma}

Here because the equation \ref{cr2} is quasi-linear, we may not get
the linear control as Lemma \ref{cont1}. We will give a self
contained proof of this result by argument similar to those of
Ahlfors and Bers. We will use their notions. In fact we will proof
the following two claims:

\begin{Claim}
\begin{equation}
d_{S^{2}}\big(h^{\mu_{n}}, h^{\mu}\big)\rightarrow 0.
\end{equation}
\end{Claim}

\begin{Proof}
Let $w^{\mu}$ be the conformal diffeomorphism described above, so we
have uniform H$\ddot{o}$lder estimates
$d_{S^{2}}\big(w^{\mu}(z_{1}), w^{\mu}(z_{2})\big)\leq c
d_{S^{2}}(z_{1}, z_{2})^{\alpha}$. Here the constant $c$ is uniform
for fixed $k<1$, when all $\|\mu\|\leq k$. Let
$h^{\mu}=(w^{\mu})^{-1}$, we have:
\begin{equation}
h^{\mu}_{\overline{w}}=\nu(w)h^{\mu}_{w},
\end{equation}

here
$\nu(w)=\big(-\mu\frac{w^{\mu}_{z}}{\overline{w^{\mu}_{z}}}\big)\circ
h$. Since $\|\nu\|_{L^{\infty}}=\|\mu\|_{L^{\infty}}$, we have
similar H$\ddot{o}$lder estimates $d_{S^{2}}\big(h^{\mu}(w_{1}),
h^{\mu}(w_{2})\big)\leq c^{\prime}
d_{S^{2}}(w_{1}, w_{2})^{\alpha}$.\\

We use contradiction arguments. Suppose $(w^{\mu_{n}})^{-1}$ do not
converge to $(w^{\mu})^{-1}$ in $L^{\infty}(S^{2}, S^{2})$, then
there exists an $\epsilon>0$ and a sequence $x_{n}\in S^{2}$ such
that $d_{S^{2}}\big((w^{\mu_{n}})^{-1}(x_{n}),
(w^{\mu})^{-1}(x_{n})\big)>\epsilon$. By the compactness of $S^{2}$,
we can assume $x_{n}\rightarrow x_{0}$, and
$(w^{\mu_{n}})^{-1}(x_{n})\rightarrow z_{1}$,
$(w^{\mu})^{-1}(x_{n})\rightarrow z_{0}$. Clearly $d_{S^{2}}(z_{0},
z_{1})\geq\epsilon$. But $w^{\mu}(z_{0})=w^{\mu}(z_{1})=x_{0}$,
which forms a contradiction since $w^{\mu}$ is a homeomorphism. This
is because of the following.\\

Denoting $z_{n}=(w^{\mu_{n}})^{-1}(x_{n})$ and
$z_{n}^{\prime}=(w^{\mu})^{-1}(x_{n})$, we have the following:
\begin{equation}
d_{S^{2}}\big(w^{\mu_{n}}(z_{n}), w^{\mu}(z_{1})\big)\leq
d_{S^{2}}\big(w^{\mu_{n}}(z_{n}),
w^{\mu_{n}}(z_{1})\big)+d_{S^{2}}\big(w^{\mu_{n}}(z_{1}),
w^{\mu}(z_{1}))\big)\rightarrow 0,
\end{equation}

The convergence of the first term is because $w^{\mu_{n}}$ have
uniform H$\ddot{o}$lder norm. So $w^{\mu}(z_{1})=x_{0}$. And
\begin{equation}
w^{\mu}(z_{0})=\lim_{n\rightarrow\infty}
w^{\mu}(z^{\prime}_{n})=\lim_{n\rightarrow\infty} x_{n}=x_{0}
\end{equation}

So we have $w^{\mu}(z_{0})=w^{\mu}(z_{1})$.
\end{Proof}

\begin{Claim}
The conformal diffeomorphism solution
$h:\mathbb{C}\rightarrow\mathbb{C}$ fixing $(0,1,\infty)$ of the
equation:
\begin{equation}
h_{w}=\alpha(w)\overline{h_{w}},
\end{equation}

have estimates:
\begin{equation}
\|(h^{\alpha}-h^{\beta})_{w}\|_{L^{p}(B_{R})}\leq
C(R)\|\alpha-\beta\|^{2\alpha}_{L^{\infty}}.
\end{equation}

Here constants depend only on bound $k$ of $|\alpha|\leq k<1$.
$\alpha=1-\frac{2}{p}$ as in Lemma \ref{cont1} and $p$ depends only
on $k$.
\end{Claim}

\begin{Proof}
We show this in five steps, and we may use $w$ to denote $h$.\\

\underline{Step 1.} We consider the following non-homogeneous
equation:
\begin{equation}
w_{\overline{z}}=\mu\overline{w_{z}}+\sigma.
\end{equation}

We want to find solutions satisfying: $w(0)=0$, $w_{z}\in
L^{p}(\mathbb{C})$, and we denote such solution by $w^{\mu,
\sigma}$. We firstly consider the following preliminary equation:
\begin{equation}
q=T(\mu\overline{q}+\sigma).
\end{equation}

Here $T$, $P$ denote the operators defined in Section 1.2 of
\cite{AB}. By the fixed point theorem, we know there is a unique
solution $q\in L^{p}(\mathbb{C})$ when $p$ is appropriate. Let
$w=P(\mu\overline{q}+\sigma)$. We have $w(0)=0$,
$w_{z}=T(\mu\overline{q}+\sigma)=q$, and
$w_{\overline{z}}=\mu\overline{q}+\sigma$ by properties of operators
$T$ and $P$ given in Lemma 3 in \cite{AB}. So
$w_{\overline{z}}=\mu\overline{w_{z}}+\sigma$, and $w$ satisfy our
restriction. So $w$ is our solution. We can know that such $w$ is
unique by estimating corresponding homogenous equation similar to
that of Lemma 1 in \cite{AB}. From properties of operators $T$ and
$P$ given in Lemma 3 in \cite{AB}, we have estimates for $w^{\mu,
\sigma}$:
\begin{equation}
\|w^{\mu, \sigma}_{z}\|_{L^{p}}=\|q\|_{L^{p}}\leq
c(p)\|\sigma\|_{L^{p}},
\end{equation}
\begin{equation}
|w^{\mu, \sigma}(z_{1})-w^{\mu, \sigma}(z_{2})|\leq
c|z_{1}-z_{2}|^{\alpha}
\end{equation}

Here $\alpha=1-\frac{2}{p}$. In fact, by the properties of $P$, we
have:
\begin{equation}
|w^{\mu, \sigma}(z_{1})-w^{\mu, \sigma}(z_{2})|\leq
c\|\mu\overline{q}+\sigma\|_{L^{p}}|z_{1}-z_{2}|^{1-\frac{2}{p}}\leq
c^{\prime}\|\sigma\|_{p}|z_{1}-z_{2}|^{1-\frac{2}{p}}.
\end{equation}

This mean our solution also have uniform H$\ddot{o}$lder norm.\\

\underline{Step 2.} $w^{\mu, \sigma}$ varies continuously in
$L^{\infty}$ and $L^{p}$ as $\mu$, $\sigma$ vary continuously. Let
$w=w^{\mu, \sigma}$, $w^{\prime}=w^{\nu, \rho}$. We have:
\begin{equation}
(w-w^{\prime})_{\overline{z}}=\mu\overline{(w-w^{\prime})_{z}}+\lambda,
\end{equation}

here $\lambda=(\mu-\nu)\overline{w^{\prime}_{z}}+(\sigma-\rho)$. By
the above results, we have estimate:
\begin{equation}
\|(w-w^{\prime})_{z}\|_{L^{p}}\leq c\|\lambda\|_{L^{p}}\leq
c\big(\|\mu-\nu\|_{L^{\infty}}+\|\sigma-\rho\|_{L^{p}}\big).
\end{equation}

Similarly, we also have estimates of H$\ddot{o}$lder norm for
$w-w^{\prime}$.\\

\underline{Step 3.} Suppose $\mu$ is compactly supported. We want to
have homeomorphism $w^{\mu}:\mathbb{C}\rightarrow\mathbb{C}$
satisfying: $w^{\mu}_{\overline{z}}=\mu\overline{w^{\mu}_{z}}$, with
normalization $w^{\mu}(0)=0$, and $w^{\mu}_{z}-1\in
L^{p}(\mathbb{C})$. In fact, let $w^{\mu}=z+w^{\mu,\mu}(z)$, with
$w^{\mu,\mu}(z)$ as in the above step, we have:
\begin{equation}
w^{\mu}_{\overline{z}}=(w^{\mu,\mu})_{\overline{z}}=\mu\overline{(w^{\mu,\mu}_{z})}+\mu=\mu(\overline{w^{\mu,\mu}_{z}+1})=\mu\overline{w^{\mu}_{z}}.
\end{equation}

Clearly, $w^{\mu}(0)=0$, and $w^{\mu}_{z}-1=w^{\mu,\mu}_{z}\in
L^{p}(\mathbb{C})$. From argument similar to Section 3.3 of
\cite{AB}, we know $w^{\mu}$ is homeomorphism. So $w^{\mu}$ is our
solution. In this case, to have a solution fixing $(0,1,\infty)$, we
only need to divide $w^{\mu}$ by $w^{\mu}(1)$. We also have results
similar to Lemma 15 in \cite{AB} that $c(R)^{-1}\leq|w^{\mu}(1)|\leq
c(R)$ when $\mu$ has compact support in $B_{R}$. We will also denote
$w^{\mu}/w^{\mu}(1)$ by $w^{\mu}$ in the following.\\

\underline{Step 4.} Let $\alpha$, $\beta$ be the two coefficients
with $|\alpha|,|\beta|\leq k<1$, we give a decomposition formula:
\begin{equation}
w^{\alpha}=w^{\beta}\circ w^{\gamma},
\end{equation}

here
$\gamma=\frac{\alpha-\beta}{1-\alpha\overline{\beta}}\frac{\beta}{\overline{\beta}}\big(\frac{(w^{\beta})^{-1}_{\overline{z}}}{(\overline{(w^{\beta})^{-1}})_{z}}\big)\circ
w^{\alpha}$. Hence $\|\gamma\|_{L^{\infty}}\leq
C\|\alpha-\beta\|_{L^{\infty}}$. The proof is just simple
calculation. Using sphere distance, we have the following estimates:
\begin{equation}
d_{S^{2}}\big(w^{\alpha}(z),
w^{\beta}(z)\big)=d_{S^{2}}\big(w^{\beta}\circ w^{\gamma}(z),
w^{\beta}(z)\big)\leq c d_{S^{2}}\big(w^{\gamma}(z),
z\big)^{\alpha}.
\end{equation}

Decompose $\gamma=\gamma_{1}+\gamma_{2}$, with $\gamma_{1}$ and
$\gamma_{2}$ supported near $0$ and $\infty$ separately, we have:
\begin{equation}
d_{S^{2}}\big(w^{\gamma}(z), z\big)\leq d_{S^{2}}\big(w^{\gamma}(z),
w^{\gamma_{1}}(z)\big)+d_{S^{2}}\big(w^{\gamma_{1}}(z), z\big).
\end{equation}

In the case $\gamma$ having compact support, for $|z|\leq R$ we
have:
\begin{equation}
d_{S^{2}}\big(w^{\gamma}(z), z\big)\leq c(R)\|w^{\gamma,
\gamma}\|_{L^{\infty}(B_{R})}=c(R)\|w^{\gamma, \gamma}-w^{\gamma,
\gamma}(0)\|_{L^{\infty}(B_{R})}\leq C(R)\|\gamma\|_{L^{\infty}}.
\end{equation}

By arguments as Section 5.1 of \cite{AB}, for $|z|\geq R$,  we have
$d_{S^{2}}\big(w^{\gamma}(z), z\big)\leq
c(R)\|\gamma\|_{L^{\infty}}$. Combining all the above together, we
have:
\begin{equation}
d_{S^{2}}\big(w^{\alpha}(z), w^{\beta}(z)\big)\leq
C(R)\|\alpha-\beta\|^{2\alpha}_{L^{\infty}}.
\end{equation}

Here sphere distance is equivalent to the ordinary plane distance
when restricted to a compact set on $\mathbb{C}$.\\

\underline{Step 5.} Choose cutoff function $\eta$ supported in
$B_{2R}$, with $\eta\equiv1$ on $B_{R}$, $\eta\leq 1$. Then we have:
\begin{equation}
\big(\eta(w^{\alpha}-w^{\beta})\big)_{\overline{z}}=\alpha\overline{\big(\eta(w^{\alpha}-w^{\beta})\big)_{z}}+\lambda,
\end{equation}

here
$\lambda=\eta(\alpha-\beta)\overline{w^{\beta}_{z}}+\eta_{\overline{z}}\big((w^{\alpha}-w^{\beta})-\alpha\overline{(w^{\alpha}-w^{\beta})}\big)$.
And $\|\lambda\|_{L^{p}}\leq
C(R)\|\alpha-\beta\|_{L^{\infty}(B_{2R})}+C^{\prime}(R)\|w^{\alpha}-w^{\beta}\|_{L^{\infty}(B_{2R})}$.
So the results in the step 1 and step 4 give:
\begin{equation}
\begin{split}
\|(w^{\alpha}-w^{\beta})_{z}\|_{L^{p}(B_{R})}&\leq\|\big(\eta(w^{\alpha}-w^{\beta})\big)_{z}\|_{L^{p}}\\
                                             &\leq C(R)\|\lambda\|_{L^{p}}\\
                                             &\leq
C(R)\|\alpha-\beta\|_{L^{\infty}(B_{2R})}+C^{\prime}(R)\|w^{\alpha}-w^{\beta}\|_{L^{\infty}(B_{2R})}.\\
                                             &\leq C(R)\|\alpha-\beta\|^{2\alpha}_{L^{\infty}}.
\end{split}
\end{equation}

Here we abuse the use of notion, and if we change $w^{\alpha}$ to
$h^{\alpha}$, and $z$ to $w$, we will get the result.
\end{Proof}

\begin{Proof}
(of Lemma \ref{cont2}) The first convergence \ref{sphere convergence
for the inverse of u-conformal map} follows from the first claim.
For the second convergence \ref{local Lp convergence of u-conformal
map}, since $h^{\mu_{n}}\rightarrow h^{\mu}$ in $L^{\infty}(S^{2},
S^{2})$, $h^{\mu_{n}}(B_{2R})$ must be restrained in a uniform
finite ball $B_{R^{\prime}}$. As $\mu_{n}\rightarrow\mu$ in
$C^{1}(\mathbb{C})$, we know
$\mu_{n}(h^{\mu_{n}}(w))\rightarrow\mu(h^{\mu}(w))$ in $L^{\infty}$
on any bounded balls $B_{2R}$ for fixed $R<\infty$. We know from the
proof of the second claim that:
\begin{equation}
\begin{split}
\|(h^{\mu_{n}}-h^{\mu})_{w}\|_{L^{p}(B_{R})}&\leq
C(R)\|\mu_{n}(h^{\mu_{n}}(w))-\mu(h^{\mu}(w))\|_{L^{\infty}(B_{2R})}\\
                                            &+C^{\prime}(R)\|h^{\mu_{n}}-h^{\mu}\|_{L^{\infty}(B_{2R})}.
\end{split}
\end{equation}

Since the sphere distance is equivalent to the plane distance on
compact sets, the first convergence result shows that
$\|h^{\mu_{n}}-h^{\mu}\|_{L^{\infty}(B_{2R})}\rightarrow 0$. So the
second convergence result \ref{local Lp convergence of u-conformal
map} holds.
\end{Proof}



Department of Mathematics, Stanford University, building 380,
Stanford, California 94305.

E-mail: xzhou08@math.stanford.edu


\end{document}